\documentclass[12pt,a4paper]{article}

\usepackage{amsmath}    
\usepackage{amssymb}
\usepackage{graphicx}   
\usepackage{verbatim}   
\usepackage{color}      

\usepackage{hyperref}   

\usepackage{enumerate}
\usepackage{mathrsfs}
\usepackage{empheq}
\usepackage{bbold}

\usepackage{multirow}%
\usepackage{amsfonts}%
\usepackage{amsthm}%
\usepackage{mathrsfs}%
\usepackage{xcolor}%
\usepackage{textcomp}%
\usepackage{manyfoot}%
\usepackage{booktabs}%
\usepackage{algorithm}%
\usepackage{algorithmicx}%
\usepackage{algpseudocode}%
\usepackage{listings}%

\usepackage[final]{showlabels}

\usepackage{amsthm}

\newtheorem{lemma}{Lemma}[section]

\newtheorem{theorem}{Theorem}[section]

\renewcommand{\theta}{\vartheta}
\renewcommand{\phi}{\varphi}

\renewcommand{\title}{A Well-Balanced Method for an Unstaggered Central Scheme, the one-space Dimensional Case}

\newcommand{\authorOne}{Yu-Chen Cheng\footnote{Institute of Mathematics, University of Wuerzburg, Wuerzburg, Germany, yu-chen.cheng@stud-mail.uni-wuerzburg.de}}
\newcommand{\authorTwo}{Christian Klingenberg\footnote{Institute of Mathematics, University of Wuerzburg, Wuerzburg, Germany, christian.klingenberg@uni-wuerzburg.de}}
\newcommand{\authorThree}{Rony Touma\footnote{Institute of Computer Science and Mathematics, Lebanese American University, Byblos, Lebanon, rony.touma@lau.edu.lb}}

\begin{document}

\begin{center} \Large
\title

\vspace{1cm}

\date{}
\normalsize

\authorOne, \authorTwo, \authorThree
\end{center}

\begin{abstract}

In this paper, we propose a new MUSCL scheme by combining the ideas of the Kurganov and Tadmor scheme and the so-called Deviation method which results in a well-balanced finite volume method for the hyperbolic balance laws, by evolving the difference between the exact solution and a given stationary solution. After that, we derive a semi-discrete scheme from this new scheme and it can be shown to be essentially TVD when applied to a scalar conservation law. In the end, we apply and validate the developed methods by numerical experiments and solve classical problems featuring Euler equations with gravitational source term.

Keywords: Euler equations, Deviation method, Unstaggered central methods, Well-balanced discretizations


\end{abstract}

\section{Introduction}
In the past few decades, central schemes for the approximate solutions of hyperbolic conservation laws were widely studied. One well-known Riemann-solver-free central scheme is the Nessyahu and Tadmor (NT) scheme introduced in \cite{ref5} which considers the approximate solutions on two staggered grids. An unstaggered version of the NT scheme is developed later in \cite{Jiang}. After that, in \cite{ref3}, Kurganov and Tadmor proposed a modified NT-type scheme (KT-scheme), which adopted narrower cells when considering the solutions over the Riemann fans; it can avoid to overestimate the values of smooth regions. Also, Authors in \cite{ref3} created a semi-discrete scheme from the previously mentioned scheme. There are some extensions based on this semi-discrete scheme: \cite{Semi1}; \cite{Semi2}; \cite{Semi3}; \cite{Semi4}.

In mechanics, the systems with outer forces have also received much attention. A system with an outer force can be written as the in-homogeneous conservation laws, which is also known as the balance laws,
\begin{equation} \label{eq1}
    \partial_t q(x,t)+\nabla_x f(q(x,t))= S(q(x,t)),
\end{equation}
where $q(x,t)=(q_1(x,t), q_2(x,t),...,q_N(x,t))^T$ is an N-vector of conserved quantities in the d-spatial variables $x=(x_1, x_2,..., x_d)$, and $f(q)=(f^1,f^2,...,f^d)$ is a nonlinear flux. $S(q)=(s^1,s^2,...,s^d)$ is a source term. There are some common body forces, for example, gravity, electric forces, magnetic forces. The so-called well-balanced schemes are the numerical methods designed to preserve discrete steady states for these kinds of systems. Well-balanced schemes can be derived from conventional numerical methods for the conservation laws by adding an appropriate discretization of the source term. There are some examples of the well-balanced scheme: \cite{ref4}; \cite{ref10}; \cite{ref14}; \cite{WB_BKL}; \cite{upwind2}; \cite{rec_FE}; \cite{rec_wb2}; \cite{relaxation3_tbc}; \cite{AMC}; \cite{WenoVol}.

In this paper, we begin by presenting a new second-order accurate unstaggered central scheme that is well-balanced by combining the concepts of the so-called Deviation method in \cite{ref14} and the KT scheme in \cite{ref3}. We firstly consider the deviation $\Delta q$ between the vector of conserved variables $q$ and a given steady state $\tilde q$, and then estimate the solutions on narrower cells over Riemann fans. This new scheme retains the advantage of being Riemann-solver-free and satisfies the well-balanced property.

On the basis of our new central scheme, we construct a corresponding semi-discrete scheme in the following section. Here we prove that this semi-discrete scheme satisfies the total-variation diminishing (TVD) property when applying to the modified homogeneous scalar conservation laws.

In the next section, we apply our fully-discrete schemes to some numerical tests related to the Euler system with gravity, which can be applied to model the physical phenomena, and show the results of the comparison between our numerical solutions and the exact solutions (or compared solution from \cite{ref4}). Finally, we end in the final section with the conclusion.

\section{A new central well-balanced strategy for in-homogeneous hyperbolic conservation laws}\label{sec:scheme}
\begin{figure}[ht] 
\centering
\includegraphics[width=0.7\textwidth]{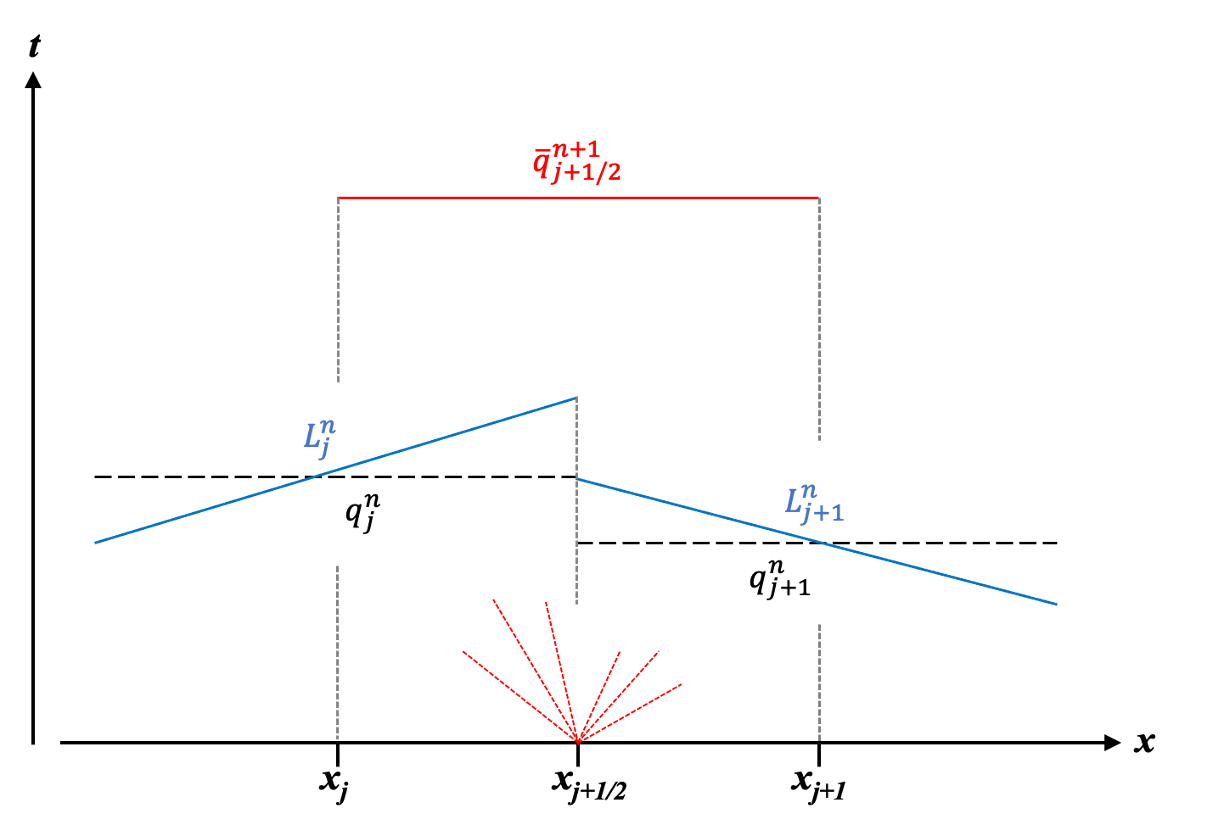}
\caption{Geometry of the NT Scheme: piecewise linear interpolants $L^n(x)$ are evolved on two staggered grids. }
\label{Fig.2.1}
\end{figure}

Kurganov and Tadmor introduced in year 2000 a modification of the Nessyahu-Tadmor scheme (NT scheme) \cite{ref5} which was originally introduced in 1990. Similarly to the NT scheme, the KT scheme evolves a piecewise linear numerical solution;  the difference between these two scheme is that in NT scheme, the fixed width control cell $[x_j,x_{j+1}]$ is considered to approximate the average over Riemann fans (see figure \ref{Fig.2.1}), while a narrower control cell is adopted in the KT scheme (see figure \ref{Fig.2.2}). The width of the control cell in KT scheme is determined by the local wave speed. Both the NT and KT schemes are suitable to approximate the solution of hyperbolic conservation laws and cannot be used alone to solve balance laws as they fail to capture steady state solutions or stationary solutions. In this work we will adopt the Deviation method from \cite{ref14} and \cite{ref4} and rewrite the Euler with gravity system in a different, yet equivalent form in terms of a target steady state solution and then we solve the resulting reformulated system with the aid of the KT scheme. We show that the resulting finite volume method preserves steady states and stationary equilibria at the discrete level while benefiting from the main characteristics and properties of the KT scheme.

\subsection{Framework of the Deviation method} \label{subsection2.2.2}
To construct our new numerical scheme for hyperbolic balance laws we firstly follow the steps in \cite{ref4} and the guidelines of the Deviation method presented in \cite{ref14} to construct the modified balanced laws, and then apply the approach of the KT scheme in section \ref{subsection2.2.1} to evolve the numerical solution.\\
Consider the one-dimensional balance laws
\begin{equation} \label{3.1}
    \left \{
    \begin{aligned}
    & q_t +f(q)_x =  S(q,x), \;\;\;  x\in \Omega \subset \mathbb{R}, t>0 \\
    & q(x,0) = q_0(x). \\
    \end{aligned}
    \right.
\end{equation}
Let $\tilde{q}$ be a given hydrostatic solution of \eqref{3.1}. In other word, it satisfies
\begin{equation} \label{3.2}
    f(\tilde{q})_x = S(\tilde{q},x).
\end{equation}
Next, subtracting \eqref{3.2} from \eqref{3.1} yields
\begin{equation} \label{3.3}
    q_t+f(q)_x-f(\tilde q)_x = S(q,x) - S(\tilde q,x).
\end{equation}
Define the deviation $\Delta q=q-\tilde{q}$. Applying $q=\Delta q+\tilde{q}$ to \eqref{3.3} leads to
\begin{equation} \label{3.4}
    (\Delta q +\tilde q)_t + [f(\Delta q+\tilde{q}) - f(\tilde{q})]_x
    = S(\Delta q+\tilde{q},x)-S(\tilde{q},x).
\end{equation}
Since $\tilde{q}$ is the stationary solution; i.e., $\frac{\partial \tilde{q}}{\partial t}=0$, it implies
\begin{equation}  \label{3.5}
    (\Delta q)_t + [f(\Delta q+\tilde{q}) - f(\tilde{q})]_x
    = S(\Delta q+\tilde{q},x)-S(\tilde{q},x).
\end{equation}
If the source term $S(q,x)$ in \eqref{3.1} is a linear functional in terms of the conserved variables,
then 
\begin{equation} \label{3.6}
    S(\Delta q+\tilde q,x) - S(\tilde q,x) = S(\Delta q,x)
\end{equation}
holds and \eqref{3.5} can be rewritten as
\begin{equation} \label{3.7}
    (\Delta q)_t + [f(\Delta q+\tilde{q}) - f(\tilde{q})]_x
    = S(\Delta q,x).
\end{equation}
This leads us to the following Lemma.

\begin{lemma}
\label{lemma3.1.1}
Consider the balance law \eqref{3.1} and a given hydrostatic solution $\tilde q$.\quad The deviation quantity $\Delta q$ satisfying the modified balance law \eqref{3.7} maintains the same local speed as the original balance law in system \eqref{3.1}.
\end{lemma}

\begin{proof}
Define
\begin{equation} \label{3.8}
    F(\Delta q) = f(\Delta q +\tilde q) -f(\tilde q).
\end{equation}
We rewrite \eqref{3.7} as
\begin{equation} \label{3.9}
    (\Delta q)_t +F(\Delta q)_x = S(\Delta q,x).
\end{equation}
Consider the flux Jacobian of \eqref{3.9}, $\frac{\partial F(\Delta q)}{\partial \Delta q}$.  According to the definition of $\Delta q$, 
\begin{equation} \label{3.10}
    \frac{\partial \Delta q}{\partial q}
    =\frac{\partial}{\partial q}  (q-\tilde q) =1.
\end{equation}
Thus, the following equality 
\begin{equation} \label{3.11}
    \frac{\partial F(\Delta q)}{\partial q} 
    = \frac{\partial F(\Delta q)}{\partial \Delta q} 
    \frac{\partial \Delta q}{\partial q}
    =\frac{\partial F(\Delta q)}{\partial \Delta q} \\
\end{equation}
holds. By the definition of $F(\Delta q)$ in \eqref{3.8},  
\begin{equation} \label{3.12}
    \frac{\partial F(\Delta q)}{\partial q}
    =\frac{\partial}{\partial q}\big(f(\Delta q +\tilde q) -f(\tilde q)\big)
    =\frac{\partial}{\partial q}f(\Delta q + \tilde q)
    =\frac{\partial f(q)}{\partial q}.
\end{equation}
holds. Hence, we obtain 
\begin{equation} \label{3.13}
    \frac{\partial F(\Delta q)}{\partial \Delta q}
    = \frac{\partial f(q)}{\partial q}.
\end{equation}
Since the local speeds are determined by the eigenvalues of the flux Jacobian, then equation \eqref{3.13} shows that the local speeds are the same as those of system \eqref{3.1}. 
\end{proof}

\subsection{Application of the Kurganov-Tadmor scheme} \label{subsection2.2.1}

Now we proceed to apply the idea of the KT scheme to the modified balance law \eqref{3.7}.  Following classical finite volume methods, we start by partitioning  the computational domain using the control cells $C_j=[x_{j-\frac{1}{2}},x_{j+\frac{1}{2}}]$ and we assume that numerical solution of the balance law \eqref{3.7} is known at time $t^n$ with $(\Delta q)^n_j $ denoting the cell-centered average value. The numerical solution at the next time step $t^{n+1}$ is obtained by the following three main steps: Reconstruction, Evolution, and Projection. Below we detail each of these steps.
\subsubsection{Reconstruction}
The finite volume method we employ in this work evolves a piecewise linear numerical solution and only the cell average values $(\Delta q)^n_j $ are stored at the centers $x_j$ of the control cells $C_j$. To avoid oscillations whenever a piecewise-linear reconstruction $\Delta q(x,t^n)$ is needed over the cells $C_j=[x_{j-\frac{1}{2}},x_{j+\frac{1}{2}}]$, we define $Q_j(x,t^n) =\Delta q(x,t^n)$ for $x\in C_j$ using
\begin{equation} \label{3.14}
    Q_j(x,t^n) = (\Delta q)^n_j + ((\Delta q)_x)^n_j(x-x_j),
\end{equation}
where $((\Delta q)_x)^n_j$ is the numerical spatial derivative. \\
The so-called $MC-\theta$ limiter is a common choice to compute the numerical derivatives in \eqref{3.14} , and it is defined as follows,
\begin{multline} \label{3.15}
    ((\Delta q)_x)^n_j 
    = \\\text{minmod} 
    \left[\theta \frac{(\Delta q)^n_{j+1}-(\Delta q)^n_j}{\Delta x}, \frac{(\Delta q)^n_{j+1}-(\Delta q)^n_{j-1}}{2\Delta x}, \theta \frac{(\Delta q)^n_j-(\Delta q)^n_{j-1}}{\Delta x} \right],
\end{multline}
with $1\leq \theta \leq2$, while the minmod function is defined as
\begin{equation} \label{3.16}
    \text{minmod}(a,b,c):= \left \{
    \begin{aligned}
        & \text{sign}(a)\min\{|a|, |b|, |c|\},\;\; \text{if}\; \text{sign}(a)=\text{sign}(b)=\text{sign}(c) \\
        & 0, \;\; \text{otherwise}.
    \end{aligned}
    \right.
\end{equation}
This $MC-\theta$ will be adopted for the numerical tests in section \ref{sec:numerics}.     

\subsubsection{Evolution}
\begin{figure}[ht] 
\centering
\includegraphics[width=1\textwidth]{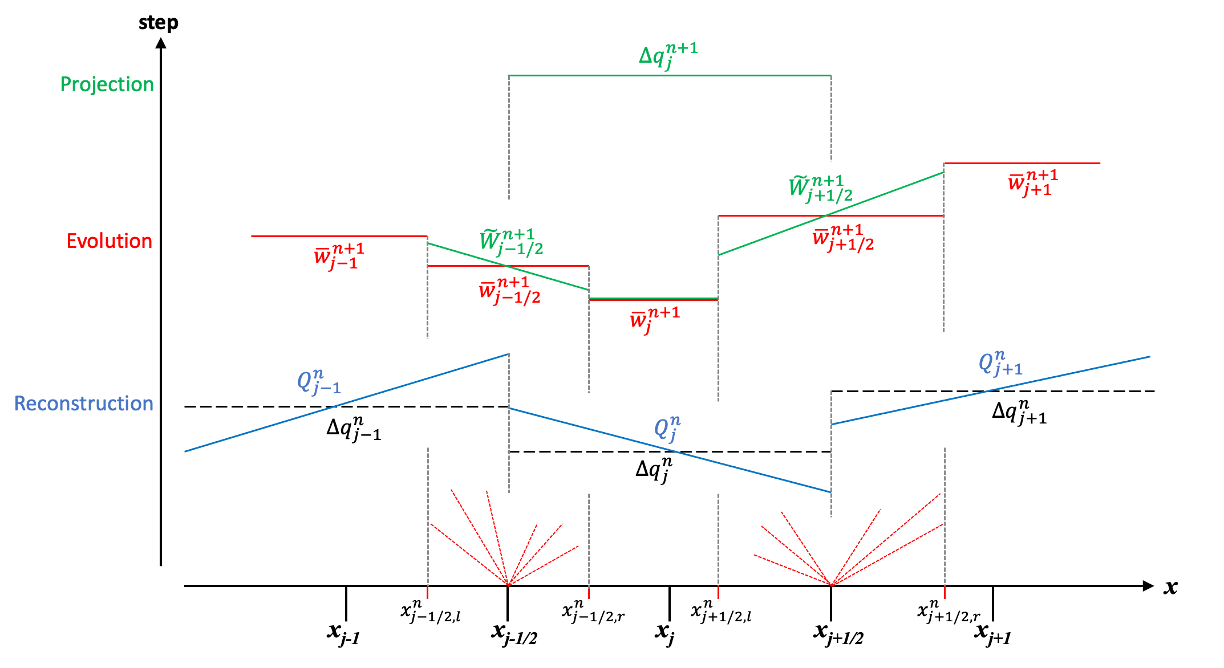}
\caption{Geometry of the KT Scheme (displayed by using the deviation $\Delta q$)}
\label{Fig.2.2}
\end{figure}

Next, we evolve the equation to the next time step $t^{n+1}=t^n+\Delta t$ and approximate the averages over the cells. Looking at the figure \ref{Fig.2.2}, the uniform cell is divided into unsmooth and smooth regions separated by the points $x^n_{j-\frac{1}{2},r}$ and $x^n_{j+\frac{1}{2},l}$, defined by

\begin{equation} 
x^n_{j+\frac{1}{2},l}:=x_{j+\frac{1}{2}}-a^n_{j+\frac{1}{2}}\Delta t \qquad \text{and} \qquad x^n_{j+\frac{1}{2},r}:=x_{j+\frac{1}{2}}+a^n_{j+\frac{1}{2}}\Delta t,
\end{equation}
where the wave speed $a^n_{j+\frac{1}{2}}$ is as defined in \cite{ref3}. In this scheme, we choose the easier way to define the wave speed: 
\begin{equation} 
    a^n_{j+\frac{1}{2}} := \max \bigg\{ \Lambda \Big( \frac{\partial f}{\partial q}(q^-_{j+\frac{1}{2}}) \Big),
    \Lambda \Big( \frac{\partial f}{\partial q}(q^+_{j+\frac{1}{2}}) \Big) \bigg\},
\end{equation}
where $q^-_{j+\frac{1}{2}}:=q^n_j+\frac{\Delta x}{2}(q_x)^n_j$ and $q^+_{j+\frac{1}{2}}:=q^n_{j+1}-\frac{\Delta x}{2}(q_x)^n_{j+1}$ are the corresponding left and right intermediate values of $q^n_{j+\frac{1}{2}}$ with the slope $(q_x)^n_j$, and $\Lambda(\frac{\partial f}{\partial q}(q))$ means the eigenvalues of the flux Jacobian $\frac{\partial f(q)}{\partial q}$.

For this reason, we consider two cases: in the first case, we discuss the average over the unsmooth interval $U^n_{j+\frac{1}{2}}=[x^n_{j+\frac{1}{2},l},x^n_{j+\frac{1}{2},r}]$; in the second case, we discuss the average over the smooth interval $M^n_j=[x^n_{j-\frac{1}{2},r},x^n_{j+\frac{1}{2},l}]$. 
\\[0.4cm]
\emph{Case 1.}\quad Unsmooth region.\\[0.2cm]
We firstly integrate \eqref{3.7} over $U^n_{j+\frac{1}{2}}\times[t^n,t^{n+1}]$, 
\begin{equation}  \label{3.18}
    \int^{t^n+1}_{t^n} \int_{U^n_{j+\frac{1}{2}}} (\Delta q)_t + [f(\Delta q+\tilde{q})-f(\tilde{q})]_x \, dxdt
    = \int^{t^n+1}_{t^n} \int_{U^n_{j+\frac{1}{2}}} S(\Delta q,x) \,dxdt.
\end{equation}
Applying the Green's theorem to the left-hand-side of \eqref{3.18} leads to
\begin{equation}  \label{3.19}
    \oint_{\partial (UT)^n_{j+\frac{1}{2}}} \left \{ [f(\Delta q+\tilde{q})-f(\tilde{q})]\; dt
    - \Delta q\; dx \right \}
    = \int^{t^{n+1}}_{t^n} \int^{x^n_{j+\frac{1}{2},r}}_{x^n_{j+\frac{1}{2},l}} S(\Delta q,x) \;dxdt, 
\end{equation}
where $(UT)^n_{j+\frac{1}{2}} = U^n_{j+\frac{1}{2}}\times [t^n,t^{n+1}]$.\\ 
Then after some calculations shown in Appendix, we obtain the intermediate value $\overline{w}^{n+1}_{j+\frac{1}{2}}$, which represents the average over the interval $U^n_{j+\frac{1}{2}}$ at $t^{n+1}$,
\begin{equation} \label{3.22}
\begin{split}
    \overline{w}^{n+1}_{j+\frac{1}{2}} 
    = \frac{(\Delta q)^n_j + (\Delta q)^n_{j+1}}{2}
    + \frac{\Delta x - \frac{\Delta x^n_{j+\frac{1}{2}}}{2}}{4} \bigg[ ((\Delta q)_x)^n_j - ((\Delta q)_x)^n_{j+1} \bigg] \\
    +  \frac{1}{\Delta x^n_{j+\frac{1}{2}}} \left [ F_U + S_U \right ],
\end{split}
\end{equation}
where   
\begin{equation} \label{3.23}
\begin{split}
    F_U
    = & - \int^{t_{n+1}}_{t_n} [f((\Delta q+\tilde{q}) (x^n_{j+\frac{1}{2},r},t))-f(\tilde{q}(x^n_{j+\frac{1}{2},r},t))]dt \\
    &+ \int^{t_{n+1}}_{t_n} [f(\Delta q+\tilde{q})(x^n_{j+\frac{1}{2},l},t))-f(\tilde{q}(x^n_{j+\frac{1}{2},l},t))]dt, 
\end{split}
\end{equation}
and
\begin{equation} \label{3.24}
    S_U = \int^{t^{n+1}}_{t^n} \int^{x^n_{j+\frac{1}{2},r}}_{x^n_{j+\frac{1}{2},l}} S(\Delta q,x) \;dxdt.
\end{equation}

Next, we consider the integral of the flux function $F_U$. Applying the midpoint rule to $F_U$, we get 
\begin{equation} \label{3.28}
\begin{split}
    F_U
    =  \Delta t & \left[ - f((\Delta q)^{n+\frac{1}{2}}_{j+\frac{1}{2},r}+\tilde q_{j+\frac{1}{2},r}) + f((\tilde q)_{j+\frac{1}{2},r}) \right. \\
      & \left. + f((\Delta q)^{n+\frac{1}{2}}_{j+\frac{1}{2},l}+\tilde q_{j+\frac{1}{2},l}) - f((\tilde q)_{j+\frac{1}{2},l}) \right]. 
\end{split}
\end{equation}
The midpoint values $(\Delta q)^{n+\frac{1}{2}}_{j+\frac{1}{2},l}$ and $(\Delta q)^{n+\frac{1}{2}}_{j+\frac{1}{2},r}$ can be approximated by the Taylor expansion in addition to the balance law:
\begin{equation} \label{3.29}
\begin{split}
    (\Delta q)^{n+\frac{1}{2}}_{j+\frac{1}{2},l}
    & = (\Delta q) ^n_{j+\frac{1}{2},l} + \frac{\Delta t}{2} ((\Delta q)_t)^n_{j+\frac{1}{2},l} \\
    & = (\Delta q) ^n_{j+\frac{1}{2},l} +\frac{\Delta t}{2} \left[ - \big[ (f(\Delta q+\tilde q)-f(\tilde q))_x \big ]_{(x^n_{j+\frac{1}{2},l},t^n)}
    + S((\Delta q)_{j+\frac{1}{2},l},t^n) \right],
\end{split}
\end{equation}
\begin{equation} \label{3.30}
\begin{split}
    & (\Delta q)^{n+\frac{1}{2}}_{j+\frac{1}{2},r} 
     = (\Delta q) ^n_{j+\frac{1}{2},r} + \frac{\Delta t}{2} ((\Delta q)_t)^n_{j+\frac{1}{2},r} \\
    &\qquad\quad = (\Delta q) ^n_{j+\frac{1}{2},r} +\frac{\Delta t}{2} \left[ - \big [ (f(\Delta q+\tilde q)-f(\tilde q))_x \big ]_{(x^n_{j+\frac{1}{2},r},t^n)}
    + S((\Delta q)_{j+\frac{1}{2},r},t^n) \right],
\end{split}
\end{equation}
and where
\begin{equation} \label{3.31}
\begin{split}
    & (\Delta q)^n_{j+\frac{1}{2},l}
    = (\Delta q)^n_j
    + \Delta x ((\Delta q)_x)^n_j \left[\frac{1}{2}-\lambda a^n_{j+\frac{1}{2}}\right], \\
    & (\Delta q)^n_{j+\frac{1}{2},r}
    = (\Delta q)^n_{j+1}
    - \Delta x ((\Delta q)_x)^n_{j+1} \left[\frac{1}{2}-\lambda a^n_{j+\frac{1}{2}}\right],
\end{split}
\end{equation}
with the mesh ratio $\lambda=\frac{\Delta t}{\Delta x}$.\\[0.1cm]
Lastly we approximate to second-order of accuracy the integral of the source term $S_U$ using the midpoint rule in time and the trapezoidal rule in space, and obtain
\begin{equation} \label{3.32}
    S_U
    = \Delta t \int^{x^n_{j+\frac{1}{2},r}}_{x^n_{j+\frac{1}{2},l}} S(\Delta q^{n+\frac{1}{2}},x) dx
    = \Delta t \Delta x^n_{j+\frac{1}{2}} \left [ \frac{S^{n+\frac{1}{2}}_{j
+\frac{1}{2},l} + S^{n+\frac{1}{2}}_{j+\frac{1}{2},r}}{2}
    \right ].
\end{equation}
Equipped with the approximation of the integral of the flux term \eqref{3.28}, and the integral of the source term \eqref{3.32}, the intermediate average over the interval $U^n_{j+\frac{1}{2}}$ reduces to
\begin{equation} \label{3.33}
\begin{split}
    \overline{w}^{n+1}_{j+\frac{1}{2}}
    =  \frac{(\Delta q)^n_j + (\Delta q)^n_{j+1}}{2}
   &  + \frac{\Delta x - \frac{\Delta x^n_{j+\frac{1}{2}}}{2}}{4} \Big[ ((\Delta q)_x)^n_j - ((\Delta q)_x)^n_{j+1} \Big] \\
    & + \frac{\Delta t}{\Delta x^n_{j+\frac{1}{2}}}  \Big[ - f((\Delta q)^{n+\frac{1}{2}}_{j+\frac{1}{2},r}+\tilde q_{j+\frac{1}{2},r}) + f((\tilde q)_{j+\frac{1}{2},r}) \\
    &\qquad\qquad\quad + f((\Delta q)^{n+\frac{1}{2}}_{j+\frac{1}{2},l}+\tilde q_{j+\frac{1}{2},l}) - f((\tilde q)_{j+\frac{1}{2},l}) \Big] \\
    & + \Delta t \Bigg [ \frac{S^{n+\frac{1}{2}}_{j
+\frac{1}{2},l} + S^{n+\frac{1}{2}}_{j+\frac{1}{2},r}}{2}
    \Bigg ].
\end{split}
\end{equation}
\emph{Case 2.} Smooth region.\\[0.1cm]
We consider the interval $M^n_j$.\quad In a similar way, integrating \eqref{3.7} over $M^n_j \times [t^n,t^{n+1}]$, we obtain
\begin{equation}  \label{3.34}
    \oint_{\partial M^n_j} \left \{ [f(\Delta q+\tilde{q})-f(\tilde{q})]\; dt
    - \Delta q\; dx \right \}
    = \int^{t^{n+1}}_{t^n} \int^{x^n_{j+\frac{1}{2},l}}_{x^n_{j-\frac{1}{2},r}} S(\Delta q,x) \,dxdt, 
\end{equation}
Define $\Delta x^n_j = x^n_{j+\frac{1}{2},l} - x^n_{j-\frac{1}{2},r}=\Delta x-\Delta t(a^n_{j-\frac{1}{2}}+a^n_{j+\frac{1}{2}})$. Then using the same calculation steps as in the appendix, the average over $M^n_j$ at $t^{n+1}$ denoted as $\overline{w}^{n+1}_j$ is approximated by  
\begin{equation}  \label{3.35}
\begin{split}
    \overline{w}^{n+1}_j 
    & := \frac{1}{\Delta x^n_j} \int^{x^n_{j+\frac{1}{2},l}}_{x^n_{j-\frac{1}{2},r}}  \Delta  q(x,t^{n+1}) dx \\
    & = \frac{1}{\Delta x^n_j} \left [
    \int^{x^n_{j+\frac{1}{2},l}}_{x^n_{j-\frac{1}{2},r}} \Delta q(x,t^n) dx + F_M +S_M \right ]   \\
    & = \frac{1}{\Delta x^n_j} 
    \left [ \Delta x^n_j Q_j(x^n_{j,m},t^n) + F_M +S_M  \right ]    \\
    & = \left[ (\Delta q)^n_j + (x^n_{j,m}-x_j) ((\Delta q)_x)^n_j  \right]
    + \frac{1}{\Delta x^n_j} \left[ F_M +S_M \right] \\
    & = (\Delta q)^n_j + \frac{\Delta x^n_{j-\frac{1}{2}}-\Delta x^n_{j+\frac{1}{2}}}{4} ((\Delta q)_x)^n_j
    + \frac{1}{\Delta x^n_j} \left[ F_M+S_M  \right],
\end{split}
\end{equation}
where $x^n_{j,m}$ denotes the midpoint of the interval $[x^n_{j-\frac{1}{2},r}, x^n_{j+\frac{1}{2},l}]$ , and where
\begin{equation}  \label{3.36}
\begin{split}
    F_M
    = & - \int^{t_{n+1}}_{t_n} [f((\Delta q+\tilde{q}) (x^n_{j+\frac{1}{2},l},t))-f(\tilde{q}(x^n_{j+\frac{1}{2},l},t))]\,dt \\
    &+ \int^{t_{n+1}}_{t_n} [f(\Delta q+\tilde{q})(x^n_{j-\frac{1}{2},r},t))-f(\tilde{q}(x^n_{j-\frac{1}{2},r},t))]\,dt, 
\end{split}
\end{equation}
and
\begin{equation} 
    S_M =\int^{t^{n+1}}_{t^n} \int^{x^n_{j+\frac{1}{2},l}}_{x^n_{j-\frac{1}{2},r}} S(\Delta q,x) \,dxdt.
\end{equation}

\noindent As in case 1, the approximations of the integral of the flux term $F_M$ and the integral of the source term $S_M$ are obtained with the help of the Taylor expansion, the balance law, the midpoint rule,  and the trapezoidal rule.  Hence, we obtain, 
\begin{equation} \label{3.38}
\begin{split}
    F_M
    =  \Delta t & \left[ - f((\Delta q) ^{n+\frac{1}{2}}_{j+\frac{1}{2},l}+\tilde q_{j+\frac{1}{2},l}) + f(\tilde q_{j+\frac{1}{2},l}) \right. \\
      & \left. + f((\Delta q) ^{n+\frac{1}{2}}_{j-\frac{1}{2},r}+\tilde q_{j-\frac{1}{2},r}) - f(\tilde q_{j-\frac{1}{2},r}) \right], 
\end{split}
\end{equation}
and
\begin{equation} \label{3.39}
    S_M
    := \Delta t \int^{x^n_{j+\frac{1}{2},l}}_{x^n_{j-\frac{1}{2},r}} S(\Delta q^{n+\frac{1}{2}},x) dx
    := \Delta t \Delta x^n_j \left [ \frac{S^{n+\frac{1}{2}}_{j-\frac{1}{2},r} + S^{n+\frac{1}{2}}_{j+\frac{1}{2},l}}{2}
    \right].
\end{equation}
Equipped with \eqref{3.35}, \eqref{3.38} and \eqref{3.39}, the intermediate average over the interval $M^n_j$ is equal to
\begin{equation} \label{3.40}
\begin{split}
    \overline{w}^{n+1}_j
    =  (\Delta q)^n_j &+ \frac{\Delta x^n_{j-\frac{1}{2}}-\Delta x^n_{j+\frac{1}{2}}}{4} ((\Delta q)_x)^n_j  \\
    & + \frac{\Delta t}{\Delta x^n_j}  \bigg[ - f((\Delta q) ^{n+\frac{1}{2}}_{j+\frac{1}{2},l}+\tilde q_{j+\frac{1}{2},l}) + f(\tilde q_{j+\frac{1}{2},l})  \\
    &\qquad\qquad + f((\Delta q) ^{n+\frac{1}{2}}_{j-\frac{1}{2},r}+\tilde q_{j-\frac{1}{2},r}) - f(\tilde q_{j-\frac{1}{2},r}) \bigg] \\
    & + \Delta t  \left [ \frac{S^{n+\frac{1}{2}}_{j-\frac{1}{2},r} + S^{n+\frac{1}{2}}_{j+\frac{1}{2},l}}{2} \right].
\end{split}
\end{equation}
\subsubsection{Projection}
Finally, the last step is to project the updated solution back over the original uniform cells $C_j=[x_{j-\frac{1}{2}},x_{j+\frac{1}{2}}]$. Due to the same reason that we want to avoid oscillations, an appropriate piecewise reconstruction is needed for the unsmooth region, and it takes the form that
\begin{equation} \label{3.41}
    \widetilde W_{j+\frac{1}{2}}(x,t^{n+1}) := \overline{w}^{n+1}_{j+\frac{1}{2}} + (x-x_{j+\frac{1}{2}})( w_x)^{n+1}_{j+\frac{1}{2}}, \qquad \forall x \in [x^n_{j+\frac{1}{2},l},x^n_{j+\frac{1}{2},r}],
\end{equation}
where the slope $w_x$ is given by
\begin{equation}
    (w_x)^{n+1}_{j+\frac{1}{2}} = \text{minmod}
    \left(
    \theta \frac{w^{n+1}_{j+\frac{1}{2}}-w^{n+1}_{j}}{x_{j+\frac{1}{2}}-x^n_{j,m}},
    \frac{w^{n+1}_{j+1}-w^{n+1}_{j}}{x^{n}_{j+1,m}-x^n_{j,m}},
    \theta \frac{w^{n+1}_{j+1}-w^{n+1}_{j+\frac{1}{2}}}{x^n_{j+1,m}-x_{j+\frac{1}{2}}}  
    \right),
\end{equation}
with $1 \leq \theta \leq 2$.
\\
Consequently, the average over $[x_{j-\frac{1}{2}},x_{j+\frac{1}{2}}]$ at time $t^{n+1}$ denoted by $(\Delta q)^{n+1}_j$ is obtained after elementary  calculations,
\begin{equation} \label{3.43}
\begin{split}
    (\Delta q)^{n+1}_j
    := \frac{1}{\Delta x} & \int^{x_{j+\frac{1}{2}}}_{x_{j-\frac{1}{2}}} \overline w \,dx  \\
    =   \frac{1}{\Delta x} & \left [ \int^{x^n_{j-\frac{1}{2},r}}_{x_{j-\frac{1}{2}}} \overline{w}^{n+1}_{j-\frac{1}{2}}  \,dx
    + \int^{x^n_{j+\frac{1}{2},l}}_{x^n_{j-\frac{1}{2},r}} \overline{w}^n_j \,dx
    + \int^{x_{j+\frac{1}{2}}}_{x^n_{j+\frac{1}{2},l}} \overline{w}^{n+1}_{j+\frac{1}{2}}\,dx    \right] \\
    =   \frac{1}{\Delta x} & \left [ \int^{x^n_{j-\frac{1}{2},r}}_{x_{j-\frac{1}{2}}} \widetilde W^{n+1}_{j-\frac{1}{2}}  \,dx
    + \int^{x^n_{j+\frac{1}{2},l}}_{x^n_{j-\frac{1}{2},r}} \overline{w}^n_j \,dx
    + \int^{x_{j+\frac{1}{2}}}_{x^n_{j+\frac{1}{2},l}} \widetilde W^{n+1}_{j+\frac{1}{2}} \,dx   \right] \\
    =  \frac{1}{\Delta x} & \left [ \frac{\Delta x^n_{j-\frac{1}{2}}}{2} 
    \left ( \overline{w}^{n+1}_{j-\frac{1}{2}} + (x^n_{j-\frac{1}{2},rm} - x_{j-\frac{1}{2}}) (w_x)^{n+1}_{j-\frac{1}{2}} \right) \right. \\
    & \qquad \qquad  +  \left ( \Delta x - \frac{\Delta x^n_{j-\frac{1}{2}} + \Delta x^n_{j+\frac{1}{2}}}{2} \right ) \overline{w}^{n+1}_j  \\
    & \left. + \frac{\Delta x^n_{j+\frac{1}{2}}}{2} 
    \left ( \overline{w}^{n+1}_{j+\frac{1}{2}} + (x^n_{j+\frac{1}{2},lm} - x_{j+\frac{1}{2}})(w_x)^{n+1}_{j+\frac{1}{2}} \right) \right]  \\
    = \frac{1}{\Delta x} & \left [ \frac{\Delta x^n_{j-\frac{1}{2}}}{2} 
    \left ( \overline{w}^{n+1}_{j-\frac{1}{2}} + \frac{\Delta x^n_{j-\frac{1}{2}}}{4} (w_x)^{n+1}_{j-\frac{1}{2}} \right) \right. \\
    & \qquad \qquad  +  \left ( \Delta x - \frac{\Delta x^n_{j-\frac{1}{2}} + \Delta x^n_{j+\frac{1}{2}}}{2} \right ) \overline{w}^{n+1}_j  \\
    & \left. + \frac{\Delta x^n_{j+\frac{1}{2}}}{2} 
    \left ( \overline{w}^{n+1}_{j+\frac{1}{2}} - \frac{\Delta x^n_{j+\frac{1}{2}}}{4} (w_x)^{n+1}_{j+\frac{1}{2}} \right) \right] \\
    = \frac{1}{\Delta x} & \left [ a^n_{j-\frac{1}{2}} \Delta t
    \left ( \overline{w}^{n+1}_{j-\frac{1}{2}} + \frac{a^n_{j-\frac{1}{2}}\Delta t}{2} (w_x)^{n+1}_{j-\frac{1}{2}} \right) \right. \\
    & \qquad \qquad  +  \left ( \Delta x - (a^n_{j-\frac{1}{2}} + a^n_{j+\frac{1}{2}})\Delta t \right ) \overline{w}^{n+1}_j  \\
    & \left. + a^n_{j+\frac{1}{2}} \Delta t 
    \left ( \overline{w}^{n+1}_{j+\frac{1}{2}} - \frac{a^n_{j+\frac{1}{2}}\Delta t}{2} (w_x)^{n+1}_{j+\frac{1}{2}} \right) \right].
\end{split}
\end{equation}

\noindent To get the desired numerical solution $q^{n+1}_j$ of the original conservation law \eqref{3.1}, we just need to add the given stationary solution $\tilde{q}_j$ as follows,
\begin{equation}
    q^{n+1}_j = (\Delta q)^{n+1}_j + \tilde q_j.
\end{equation}
Next, we demonstrate that the new scheme we presented satisfies the well-balanced property, which means that the computed numerical solution $q^n_j$ remains identical to the stationary solution $ \tilde{q}_j$, $\forall j$, whenever the initial condition corresponds to the stationary solution ($q^0_j=\tilde{q}_j$ ).\\
By the definition of the deviation quantity $(\Delta q)^n_j$, if the initial condition is such that $q^0_j=\tilde{q}_j$, then this implies that $(\Delta q)^0_j=0$. Subsequently, both $(\Delta q)^0_{j+\frac{1}{2},r}$ and $(\Delta q)^0_{j+\frac{1}{2},l}$ in \eqref{3.31} are equal to zero. This leads to that $(\Delta q)^{0+\frac{1}{2}}_{j+\frac{1}{2},r}=0$ and $(\Delta q)^{0+\frac{1}{2}}_{j+\frac{1}{2},l}=0$.\\  
Then the flux approximations in \eqref{3.28} and \eqref{3.38} reduce to
\begin{equation}
    F^0_U = \Delta t \left  [ -f(\tilde q_{j+\frac{1}{2},r}) + f(\tilde q_{j+\frac{1}{2},r})
    + f(\tilde q_{j+\frac{1}{2},l}) - f(\tilde q_{j+\frac{1}{2},l})  \right ]
    =0,
\end{equation}
\begin{equation}
    F^0_M = \Delta t \left  [ -f(\tilde q_{j+\frac{1}{2},l}) + f(\tilde q_{j+\frac{1}{2},l})
    + f(\tilde q_{j-\frac{1}{2},r}) - f(\tilde q_{j-\frac{1}{2},r})  \right ]
    =0,
\end{equation}
and the estimated values of the source term $S^0_U$ and $S^0_M$ both vanish as well.  Hence, the intermediate averages in \eqref{3.33} and \eqref{3.40} reduce to
\begin{equation}
    \Delta w^{0+1}_{j+\frac{1}{2}} =0,
\end{equation}
\begin{equation}
    \Delta w^{0+1}_j =0.
\end{equation}
Then the deviation $(\Delta q)_j^1$ is identically zero at time $t=t^1$.
Following a similar reasoning, we can show that $(\Delta q)^2_j=0$, $(\Delta q)^3_j=0$, and so on. Therefore, we can deduce that $(\Delta q)_j^n$ is equal to zero at time $t=t^n$ of any subsequent iteration. \\ Finally, since $(\Delta q)_j^n=0, \forall t=t^n$, this implies
\begin{equation}
    q_j^n = (\Delta q)_j^n+\tilde{q}_j = \tilde{q}_j, \quad \forall t=t^n.
\end{equation}
Thus, we conclude that the computed numerical solution $q_j^n$ remains stationary and we proved that our scheme satisfies the well-balanced property.

\section{Semi-discrete scheme}\label{sec:semi}
In this section, we construct a semi-discrete scheme from the fully-discrete scheme \eqref{3.43}, and then show the TVD property of the semi-discrete scheme. For the proof we use a homogeneous scalar conservation law, which we discretize using the Deviation method. 
\par
Inspired by \cite{ref3}, in order to construct a semi-discrete scheme, we firstly compute the value of $\frac{(\Delta
 q)^{n+1}_j-(\Delta q)^n_j}{\Delta t}$, and then let $\Delta t \to 0$, i.e.,
\begin{equation}
    \frac{d}{dt}(\Delta q)_j(t)
    = \lim_{\Delta t \to 0} \frac{(\Delta q)^{n+1}_j-(\Delta q)^n_j}{\Delta t}.
\end{equation}
Substituting the result in \eqref{3.43} for $(\Delta q)^{n+1}_j$ (the term $\mathcal{O}(\Delta t)$ denotes all terms that are proportional to $\Delta t$) yields   
\begin{equation}\label{3.52}
\begin{split}
    & \frac{(\Delta q)^{n+1}_j-(\Delta q)^n_j}{\Delta t} =  \\
    & = \frac{a^n_{j-\frac{1}{2}}}{\Delta x} \overline{w}^{n+1}_{j-\frac{1}{2}}
    + \left(\frac{1}{\Delta t}-\frac{a^n_{j-\frac{1}{2}}+a^n_{j+\frac{1}{2}}}{\Delta x}\right) \overline{w}^{n+1}_j
    +\frac{a^n_{j+\frac{1}{2}}}{\Delta x} \overline{w}^{n+1}_{j+\frac{1}{2}}
    -\frac{1}{\Delta t} (\Delta q)^n_j + \mathcal{O}(\Delta t)   \\
    & = - \frac{1}{2\Delta x}  \left[ f(q^{n+\frac{1}{2}}_{j+\frac{1}{2},r}+\tilde q_{j+\frac{1}{2},r}) - f(\tilde q_{j+\frac{1}{2},r})
    + f(q^{n+\frac{1}{2}}_{j+\frac{1}{2},l}+\tilde q_{j+\frac{1}{2},l}) - f(\tilde q_{j+\frac{1}{2},l})  \right.  \\
    & \qquad\qquad\qquad \left. - f(q^{n+\frac{1}{2}}_{j-\frac{1}{2},r}+\tilde q_{j-\frac{1}{2},r}) + f(\tilde q_{j-\frac{1}{2},r})
    - f(q^{n+\frac{1}{2}}_{j-\frac{1}{2},l}+\tilde q_{j-\frac{1}{2},l}) + f(\tilde q_{j-\frac{1}{2},l}) \right]  \\
    &\quad + \frac{a^n_{j+\frac{1}{2}}}{2\Delta x} \left[ ((\Delta q)^n_{j+1}-\frac{\Delta x}{2}((\Delta q)_x)^n_{j+1}) - ((\Delta q)^n_j+\frac{\Delta x}{2}((\Delta q)_x)^n_j) \right]  \\
    &\quad - \frac{a^n_{j-\frac{1}{2}}}{2\Delta x} \left[ ((\Delta q)^n_j-\frac{\Delta x}{2}((\Delta q)_x)^n_j) - ((\Delta q)^n_{j-1}+\frac{\Delta x}{2}((\Delta q)_x)^n_{j-1}) \right]  \\
    &\quad + \frac{S^{n+\frac{1}{2}}_{j+\frac{1}{2},l}-S^{n+\frac{1}{2}}_{j-\frac{1}{2},r}}{2}  + \mathcal{O}(\Delta t).
\end{split}
\end{equation}
By the definition of the midpoint values in \eqref{3.29}, \eqref{3.30} and \eqref{3.31}, we have
\begin{equation}
\begin{split}
    & (\Delta q)^{n+\frac{1}{2}}_{j+\frac{1}{2},r} \\
    & = (\Delta q)^n_{j+1} - \Delta x ((\Delta q)_x)^n_{j+1}(\frac{1}{2}-\lambda a^n_{j+\frac{1}{2}}) 
    - \frac{\Delta t}{2} [f((\Delta q)^n_{j+\frac{1}{2},r}+\tilde  q_{j+\frac{1}{2},r})-f(\tilde q_{j+\frac{1}{2},r})] \\
    & (\Delta q)^{n+\frac{1}{2}}_{j+\frac{1}{2},l} \\
    & = (\Delta q)^n_j + \Delta x ((\Delta q)_x)^n_j(\frac{1}{2}-\lambda a^n_{j+\frac{1}{2}}) 
    - \frac{\Delta t}{2} [f((\Delta q)^n_{j+\frac{1}{2},l}+\tilde q_{j+\frac{1}{2},l})-f(\tilde q_{j+\frac{1}{2},l})]. 
\end{split}
\end{equation}
Applying the Taylor expansion to $\tilde q_{j+\frac{1}{2},r}$ and $\tilde q_{j+\frac{1}{2},l}$, we obtain
\begin{equation}
    \tilde q_{j+\frac{1}{2},r}:=\tilde q_{j+\frac{1}{2}} + a^n_{j+\frac{1}{2}}\Delta t (\tilde q_x)_{j+\frac{1}{2}}, 
    \quad\text{and}\quad
    \tilde q_{j+\frac{1}{2},l}:=\tilde q_{j+\frac{1}{2}} - a^n_{j+\frac{1}{2}}\Delta t (\tilde q_x)_{j+\frac{1}{2}}.
\end{equation}
As $\Delta t\to 0$, the midvalues approach
\begin{equation}
\begin{split}
    (\Delta q)^{n+\frac{1}{2}}_{j+\frac{1}{2},r}
    & \to (\Delta q)_{j+1}(t) - \frac{\Delta x}{2}((\Delta q)_x)^n_{j+1} 
    =: (\Delta q)^+_{j+\frac{1}{2}}(t), \\
     (\Delta q)^{n+\frac{1}{2}}_{j+\frac{1}{2},l}
    & \to (\Delta q)_j(t) + \frac{\Delta x}{2}((\Delta q)_x)^n_j 
    =: (\Delta q)^-_{j+\frac{1}{2}}(t), \\
     \tilde q_{j+\frac{1}{2},r} & \to \tilde q_{j+\frac{1}{2}},\\
     \tilde q_{j+\frac{1}{2},l}&  \to \tilde q_{j+\frac{1}{2}}.
\end{split}
\end{equation}
Hence, letting $\Delta t \to 0$ in \eqref{3.52}, the semi-discrete scheme takes form as
\begin{equation} \label{3.56}
\begin{split}
     \frac{d}{dt}(\Delta q)_j(t)
    & = \lim_{\Delta t \to 0} \frac{(\Delta q)^{n+1}_j-(\Delta q)^n_j}{\Delta t} \\
    & = -\frac{1}{2\Delta x} \bigg[ \Big( f((\Delta q)^-_{j+\frac{1}{2}}(t)+\tilde q_{j+\frac{1}{2}}(t)) - f(\tilde q_{j+\frac{1}{2}}(t)) \\
    &\qquad\qquad\qquad\qquad\qquad + f((\Delta q)^+_{j+\frac{1}{2}}(t)+\tilde q_{j+\frac{1}{2}}(t)) - f(\tilde q_{j+\frac{1}{2}}(t)) \Big)    \\
    &\qquad\quad -  \Big( f((\Delta q)^-_{j-\frac{1}{2}}(t)+\tilde q_{j-\frac{1}{2}}(t)) - f(\tilde q_{j-\frac{1}{2}}(t)) \\
    &\qquad\qquad\qquad\qquad\qquad + f((\Delta q)^+_{j-\frac{1}{2}}(t)+\tilde q_{j-\frac{1}{2}}(t)) - f(\tilde q_{j-\frac{1}{2}}(t)) \Big)    \bigg] \\
    & \quad +\frac{1}{2\Delta x} \bigg[ a^n_{j+\frac{1}{2}} ((\Delta q)^+_{j+\frac{1}{2}}(t)-(\Delta q)^-_{j+\frac{1}{2}}(t)) \\
    & \qquad\qquad\qquad\qquad\qquad\qquad\quad -a^n_{j-\frac{1}{2}}((\Delta q)^+_{j-\frac{1}{2}}(t)-(\Delta q)^-_{j-\frac{1}{2}}(t))    \bigg ] \\
    & \quad + \frac{S((\Delta q)^-_{j+\frac{1}{2}})+S((\Delta q)^+_{j-\frac{1}{2}})}{2}.
\end{split}
\end{equation}
Then we reformulate \eqref{3.56} as
\begin{equation} \label{3.57}
    \frac{d}{dt}(\Delta q)_j(t)=-\frac{H_{j+\frac{1}{2}}(t)-H_{j-\frac{1}{2}}(t)}{\Delta x}
    + S_j(t),
\end{equation}
with the numerical flux
\begin{equation} \label{3.58}
\begin{split}
    H_{j+\frac{1}{2}}&(t)  := \\
    & \frac{1}{2} \left[ F((\Delta q)^-_{j+\frac{1}{2}})(t) + F((\Delta q)^+_{j+\frac{1}{2}})(t) \right] 
    - \frac{a_{j+\frac{1}{2}}(t)}{2} \left( (\Delta q)^+_{j+\frac{1}{2}}(t)-(\Delta q)^-_{j+\frac{1}{2}}(t) \right),
\end{split}
\end{equation}
where
\begin{equation} \label{3.59}
    F((\Delta q)^\mp_{j+\frac{1}{2}}) 
    := f((\Delta q)^\mp_{j+\frac{1}{2}}+\tilde q_{j+\frac{1}{2}}) - f(\tilde q_{j+\frac{1}{2}}),
\end{equation}
and the source term
\begin{equation}
    S_j(t) := \frac{S((\Delta q)^-_{j+\frac{1}{2}})+S((\Delta q)^+_{j-\frac{1}{2}})}{2}.
\end{equation}

\begin{theorem}(TVD of semi-discrete scheme for the modified homogeneous conservation laws)
    Consider the semi-discrete scheme \eqref{3.57} in conservation form
    \begin{equation} \label{3.57.c}
        \frac{d}{dt}(\Delta q)_j(t)=-\frac{H_{j+\frac{1}{2}}(t)-H_{j-\frac{1}{2}}(t)}{\Delta x}.
    \end{equation}
    Assume that the numerical spatial derivatives $(\Delta q_x)_j(t)$ are chosen as \eqref{3.15}.
Then the scheme \eqref{3.57.c} satisfies the TVD property.
\end{theorem}

\begin{proof}
    Consult the proof of Theorem 4.1 in \cite{ref3}, and example 2.4 in \cite{ref6}. 
\end{proof}

It has been shown that $\Delta q$ is TVD. With this property, we can broadly ensure that the full solution $q$ ($=\Delta q+\tilde q$) is non-oscillatory.  
\\
\\
\textbf{Remark.}
In this remark, we discuss that if a given finite volume method satisfies the TVD property, to what extent the numerical method combining this finite volume method with the Deviation method will also satisfy the TVD property.

In general, for a given stationary solution, the finite volume residual is
\begin{equation} \label{eq57}
    Res(\tilde q_j) = -\frac{1}{\Delta x}[
    F(\tilde q^L_{j+\frac{1}{2}},\tilde q^R_{j+\frac{1}{2}})
    -F(\tilde q^L_{j-\frac{1}{2}},\tilde q^R_{j-\frac{1}{2}})
    ] 
    =\mathcal{O}((\Delta x)^N),
\end{equation}
which expresses the local truncation error of the finite volume method evaluated on the stationary solution. Here $N$ denotes the order of accuracy of the method.

Consider the Deviation method with forward Euler method in time,
\begin{equation} \label{eq58}
    q^{n+1}_j 
    = q^n_j
    - \frac{\Delta t}{\Delta x}[
    F(q^L_{j+\frac{1}{2}},  q^R_{j+\frac{1}{2}})
    - F(q^L_{j-\frac{1}{2}},  q^R_{j-\frac{1}{2}})
    ]   -\Delta t Res(\tilde q_j),
\end{equation}
where
\begin{equation}
    q^{L(R)}_{j+\frac{1}{2}} = \Delta q^{L(R)}_{j+\frac{1}{2}} +\tilde q^{L(R)}_{j+\frac{1}{2}}.
\end{equation}
Set $\hat{Q}$ is the solution of a FV method as follows
\begin{equation} \label{eq60}
    \hat{Q}^{n+1}_j 
    = q^n_j
    - \frac{\Delta t}{\Delta x}[
    F(q^L_{j+\frac{1}{2}},  q^R_{j+\frac{1}{2}})
    - F(q^L_{j-\frac{1}{2}},  q^R_{j-\frac{1}{2}})
    ].
\end{equation}
Then comparing \eqref{eq58} and \eqref{eq60}, we obtain the relation
\begin{equation} \label{eq61}
    q^{n+1}_j = \hat{Q}^{n+1}_j -\Delta tRes(\tilde q_j).
\end{equation}
Assume $\hat{Q}$ satisfies the TVD property; i.e.,
\begin{equation} \label{eq62}
    TV(\hat{Q}^{n+1}) \leq TV(\hat{Q}^n).
\end{equation}
By equation \eqref{eq61}, we have
\begin{equation} \label{eq63}
\begin{split}
    TV(q^{n+1}) -TV(q^n)
    & = \sum_j |q^{n+1}_{j+1} - q^{n+1}_j| - \sum_j |q^n_{j+1} - q^n_j| \\
    & = \sum_j |(\hat{Q}^{n+1}_{j+1} - \Delta tRes(\tilde q_{j+1})) - (\hat{Q}^{n+1}_j - \Delta tRes(\tilde q_j))| \\
    &\qquad - \sum_j |(\hat{Q}^n_{j+1} - \Delta tRes(\tilde q_{j+1})) - (\hat{Q}^n_j - \Delta tRes(\tilde q_j))| \\
    & = \sum_j |\hat{Q}^{n+1}_{j+1}  - \hat{Q}^{n+1}_j - \Delta t(Res(\tilde q_{j+1}) - Res(\tilde q_j))| \\
    &\qquad - \sum_j |\hat{Q}^n_{j+1}  - \hat{Q}^n_j - \Delta t(Res(\tilde q)_{j+1} - Res(\tilde q_j))| \\
    &\leq \underbrace{ \sum_j |\hat{Q}^{n+1}_{j+1}  - \hat{Q}^{n+1}_j| - \sum_j |\hat{Q}^n_{j+1}  - \hat{Q}^n_j| }_{\leq0, by \eqref{eq62}}
    + \mathcal{O}(\Delta t \Delta x^N).
\end{split}
\end{equation}
For a small grid size, the error term $\mathcal{O}(\Delta t \Delta x^N)$ is quite small. This implies that in \eqref{eq63} the solution $q$ satisfies the TVD property up to at most an order one higher than the order of the underlying scheme. So it is essentially TVD. This is in the spirit similar to \cite{ENO}.

In conclusion, since Theorem 4.1 in \cite{ref3} has shown that the KT semi-discrete scheme has the TVD property, we can implies that our semi-discrete scheme (essentially) has this non-oscillatory property.    

\section{Numerical experiments and validations}\label{sec:numerics}
In this section, we present the results of a number of numerical experiments in the Euler system with gravitational source term by employing our designed well-balanced scheme in section \ref{sec:scheme}.\\
The one-dimensional Euler system with gravitational source term is given by
\begin{equation} 
    \left \{
    \begin{aligned}
    & q_t +f(q)_x =  S(q,x), \;\;\;  x\in \Omega \subset \mathbb{R}, t>0 \\
    & q(x,0) = q_0(x), \\
    \end{aligned}
    \right.
\end{equation}
where
\begin{equation} 
    q =
    \begin{pmatrix}
        \rho \\ \rho u \\ \; E \;
    \end{pmatrix},\quad
    f(q) =
    \begin{pmatrix}
        \rho u \\ \rho u^2+p \\ \; (E+p)u \;
    \end{pmatrix},\quad
    S(u)=
    \begin{pmatrix}\quad
        0 \\ -\rho\phi_x  \\ \; -\rho u \phi_x \;
    \end{pmatrix}.
\end{equation}
$\rho$, $u$, $\rho u$, and $p$ are used to denote the fluid’s density, velocity, momentum, and pressure, respectively. E is the total energy. The given function $\phi=\phi(x)$ is the gravitational field and the ratio of the specific heats $\gamma$ is suggested to be 1.4 for an ideal gas. \\
In all the tests, the parameter of the $(MC-\theta)$ limiter  is set as $\theta=1.5$. With the help of the CFL condition, the computation of the time-step is defined by
\begin{equation}
    \Delta t  = \text{CFL}\cdot \frac{\Delta x}{\max(|\lambda_k|)}.
\end{equation}
Here, the $\lambda_k$'s denote the eigenvalues of the flux Jacobian $\frac{\partial f(q)}{\partial q}$ and the CFL number is taken to be 0.485 for all test cases.\\
For some of the numerical experiments we consider here, we compare our solutions with the solutions from \cite{ref4}. The scheme considered in \cite{ref4} is the combination of the NT scheme and the Deviation method.

\subsection{Isothermal Equilibrium}\label{subsection5.1.1}
In the first numerical experiment, we consider the steady isothermal state with a linear gravitational field $\phi_x=g=1$ (see \cite{ref10}). The isothermal equilibrium state is given by
\begin{equation}
\begin{split}
    \rho(x) &= \rho_0 \text{exp}(-\frac{\rho_0g}{p_0}x),\\
    u(x) &= 0,\\
    p(x) &= p_0\text{exp}(-\frac{\rho_0g}{p_0}x).
\end{split}
\end{equation}
We set $\rho_0=1$ and $p_0=1$. The chosen stationary solution $\tilde q$ is the isothermal equilibrium state. The solution is computed on  200 grid points of the interval $[0,1]$ with the outflow boundary condition until the final time $t=0.25$ and we compared it with the exact solution.
\begin{figure}[th]
  \centering
  \begin{tabular}{ c @{\quad} c }
    \includegraphics[width=18em]{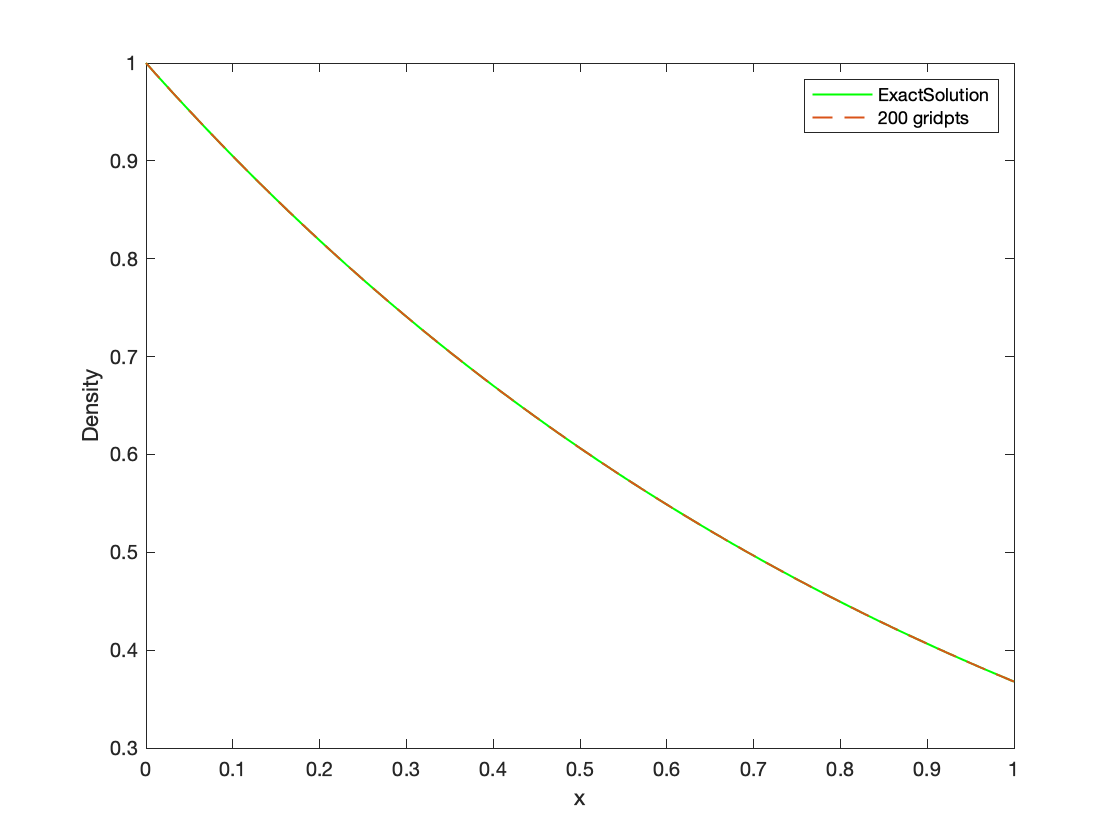} &
      \includegraphics[width=18em]{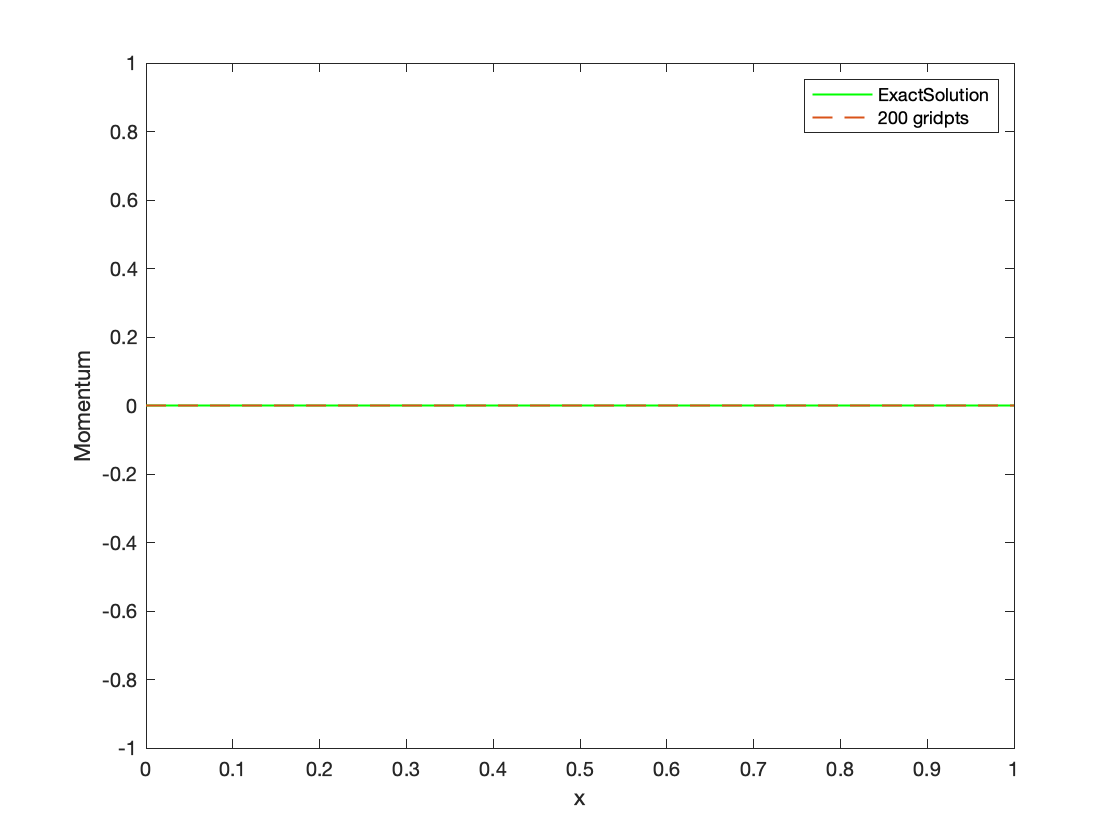} \\
    \small Density &
      \small Momentum\\
     \includegraphics[width=18em]{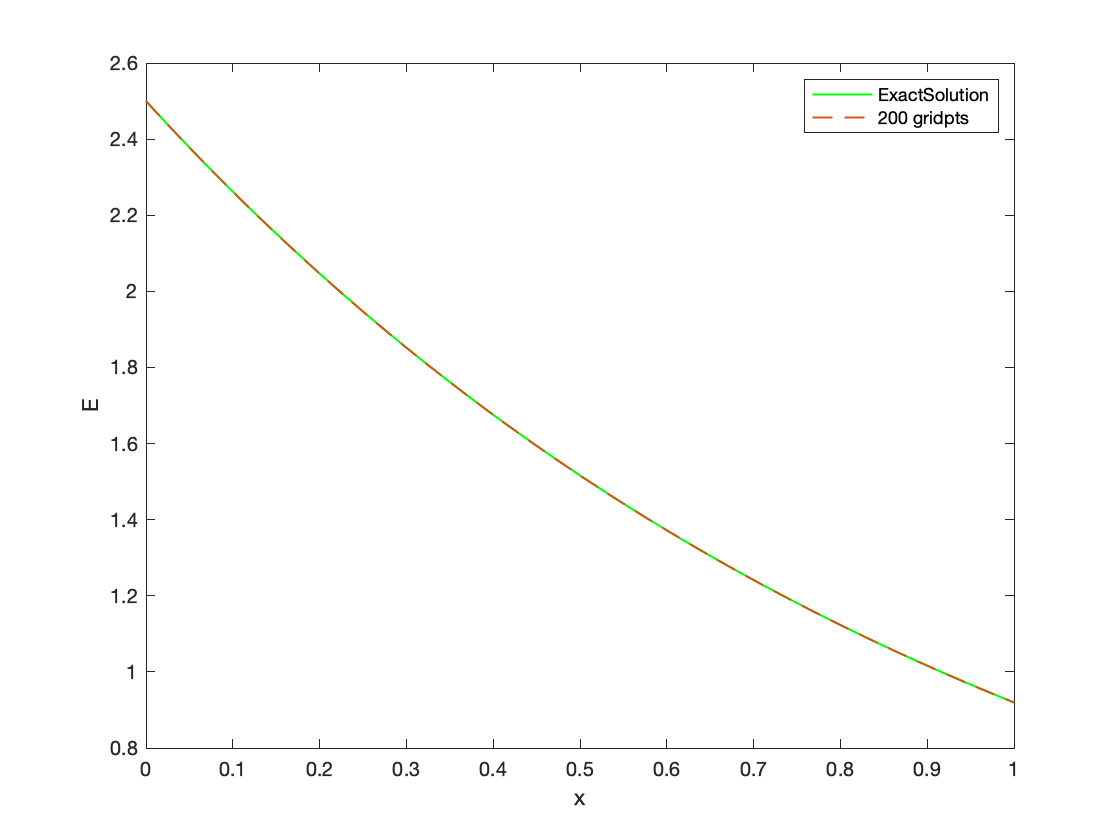} &
      \includegraphics[width=18em]{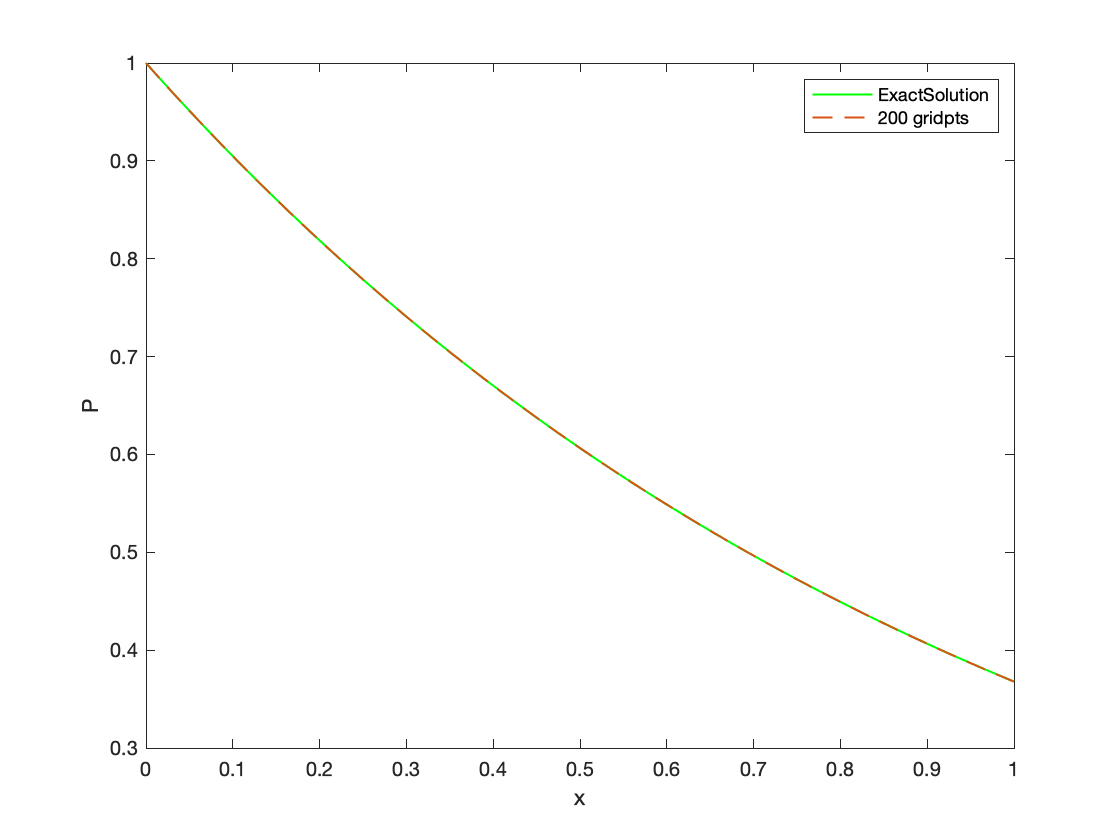} \\
    \small Energy&
      \small Pressure 
  \end{tabular}
  \vspace*{8pt}
  \caption{Results of 1D isothermal equilibrium.\label{fig.5.1-1}}
\end{figure}
 Figure \ref{fig.5.1-1} depicts the obtained results of the density, momentum, energy, and pressure which we compared with the exact equilibrium solution. The depicted results show a perfect match, thus validating the well-balanced property of the constructed finite volume method. 
\subsection{Isothermal Equilibrium with Perturbation} \label{subsection5.1.2}
\begin{figure} [th]
    \centering
    \includegraphics[width=0.6\textwidth]{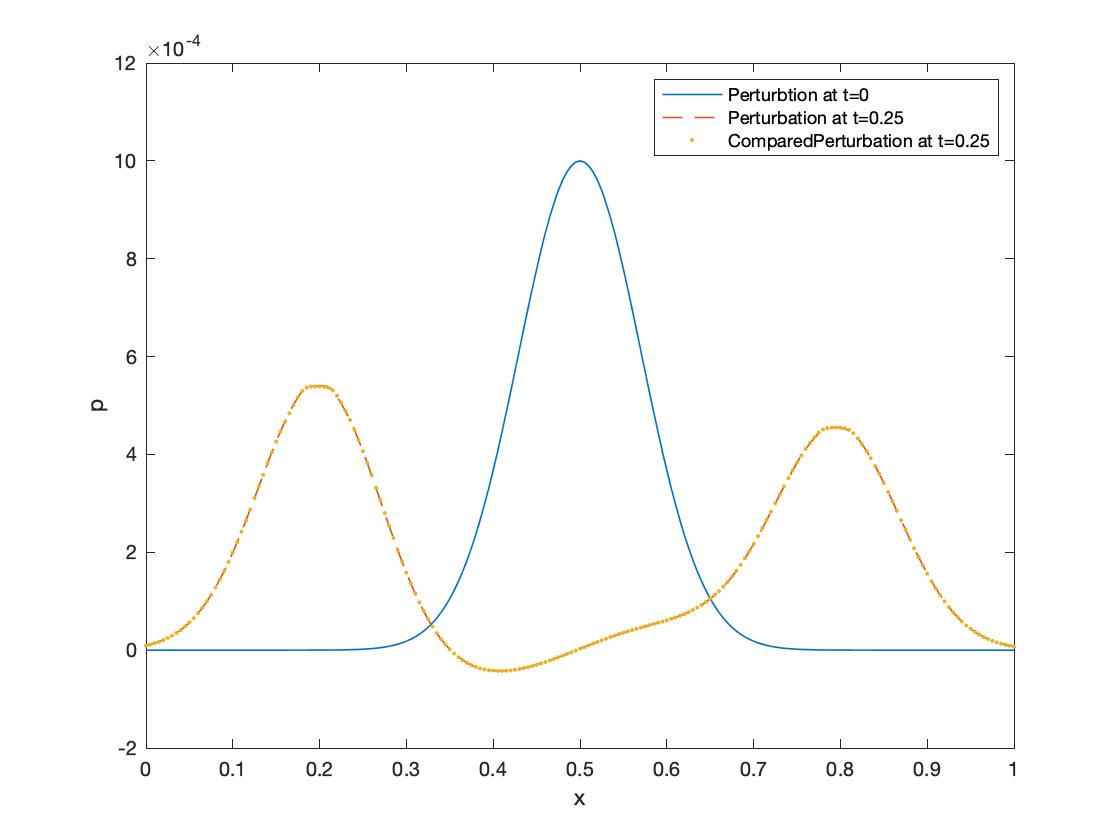}
    \caption{Initial perturbation at t=0 compared to the perturbation at the final time t=0.25 and the perturbation from compared solution. }
    \label{fig.5.1-2}
\end{figure}

The second numerical experiment is a modification to the first one, which consists of adding a small perturbation to the initial pressure given by
\begin{equation}
    p(x) = p_0\text{exp}(-\frac{\rho_0g}{p_0}x) + \eta \text{exp}(-100\frac{\rho_0g}{p_0}(x-0.5)^2),
\end{equation}
with $\eta =0.001$. Figure \ref{fig.5.1-2} shows the comparison of the initial perturbation at $t=0$, with the final perturbation at $t=0.25$ on 200 grid points obtained using the proposed scheme, and also the final perturbation at $t=0.25$ on 200 grid points obtained using the numerical scheme developed in \cite{ref4}. The obtained results show an excellent match. 
In this test case  we use $L_1$-norm to compute the errors in the density, pressure and the energy, and report these  errors and the corresponding convergence rates in table \ref{tab.5.1}. 
The reported results in table \ref{tab.5.1} confirm the second-order of accuracy of the proposed numerical scheme. 
\begin{table}[th]
    \centering
    \begin{tabular}[c]{|c|c|c|c|c|c|c|}
    \hline
    N & $\rho$ $L_1$-error & rate & $p$ $L_1$-error &  rate & $E$ $L_1$-error &  rate \\
    \hline
    200 & 3.3030E-06& - & 4.4358E-06 & - & 1.1091E-05 & - \\
    \hline
    400 & 1.4317E-06 & 1.21 & 1.9702E-06 & 1.17& 4.9260E-0.6 & 1.17\\
    \hline
    800 & 5.2586E-07 & 1.44 & 7.3033E-07 & 1.43& 1.8260E-06 & 1.43\\
    \hline
    1600 & 8.4609E-08 & 2.64 & 1.1739E-07 & 2.64 & 2.9351E-07 & 2.64 \\
    \hline
    \end{tabular}
    \caption{1D isothermal equilibrium with perturbation: $L_1$-errors and convergence rates.}
    \label{tab.5.1}
\end{table}

\subsection{Moving Equilibrium}\label{subsection5.1.3}
\begin{figure}[th]
  \centering
  \begin{tabular}{ c @{\quad} c }
    \includegraphics[width=17em]{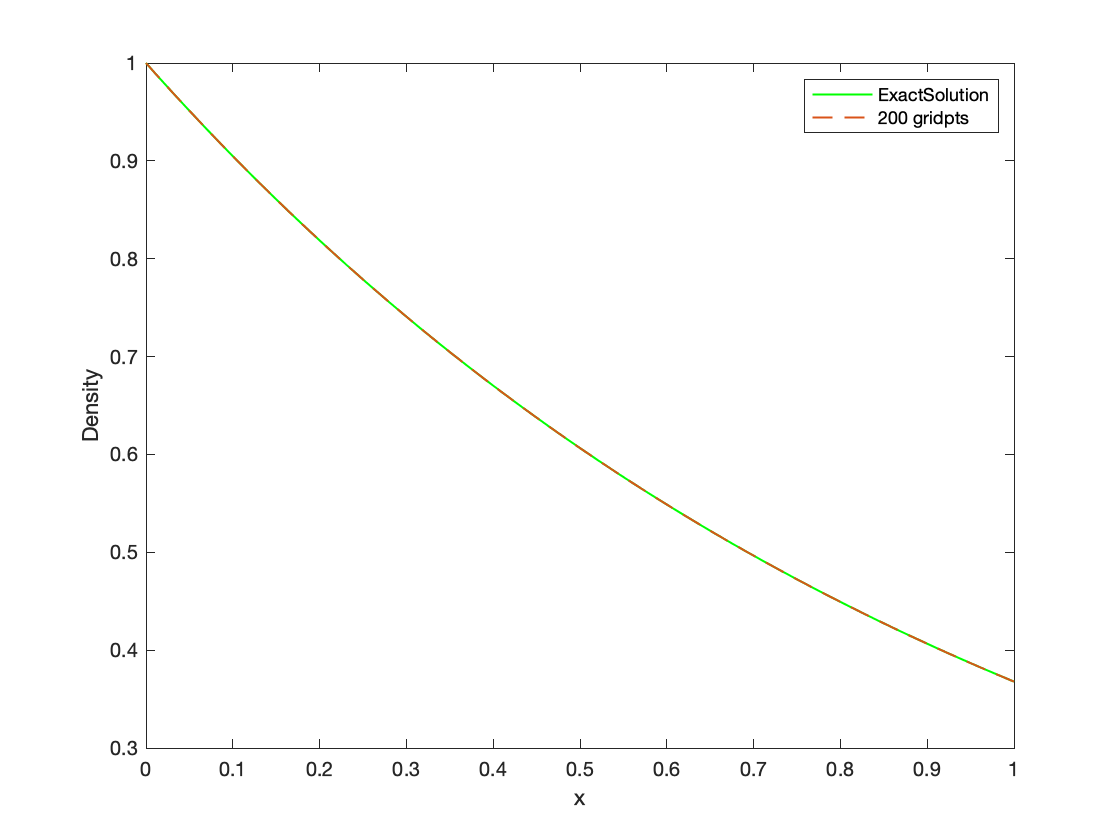} &
      \includegraphics[width=17em]{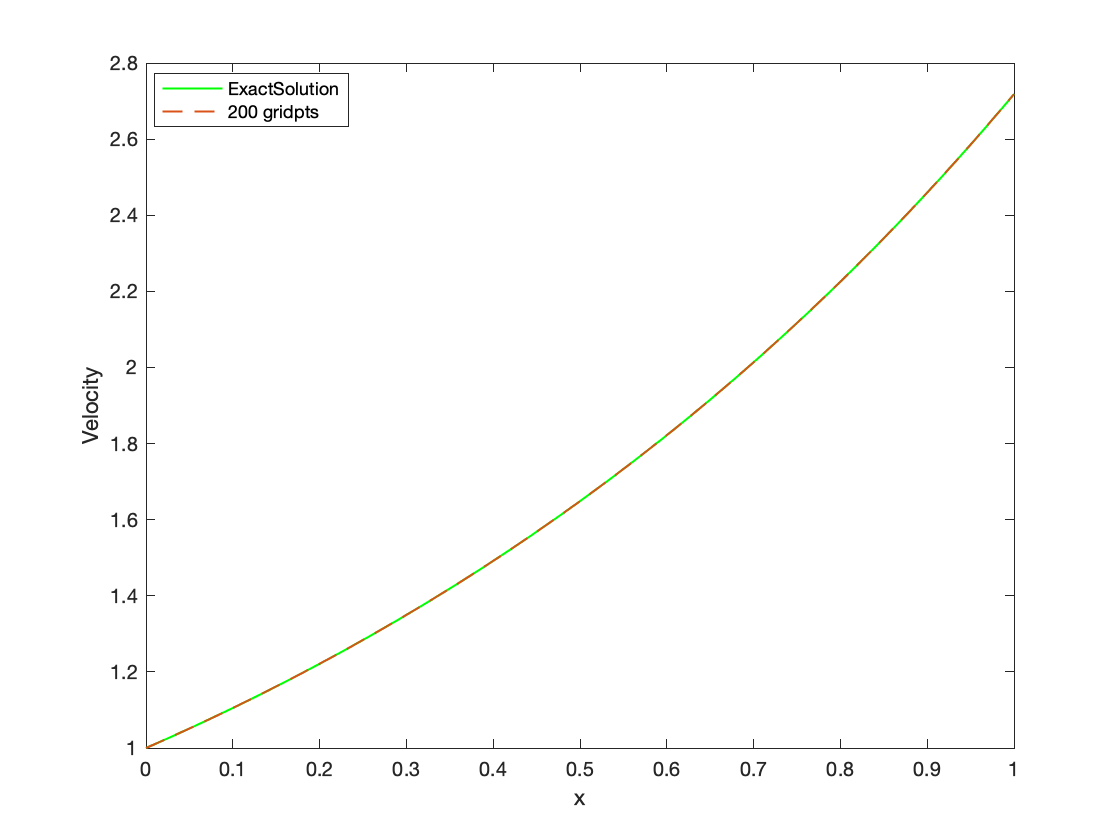} \\
    \small Moving Equilibrium: Density &
      \small Moving Equilibrium: Velocity \\
     \includegraphics[width=17em]{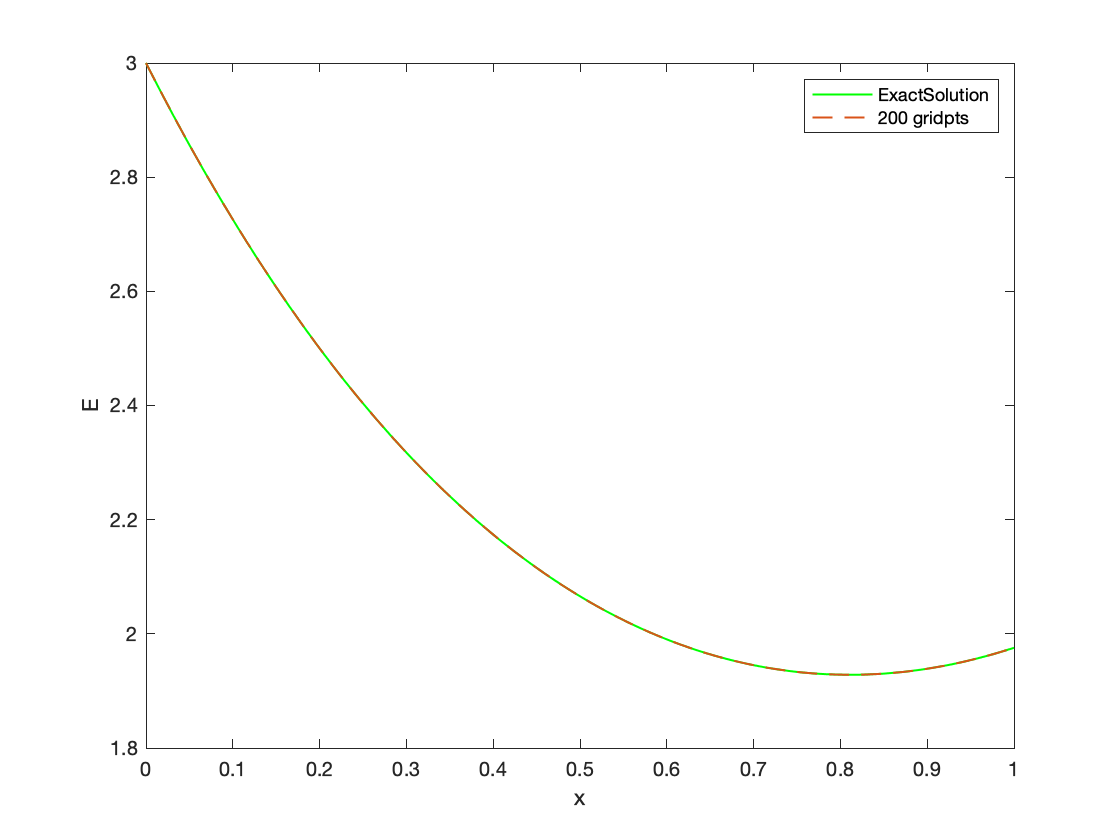} &
      \includegraphics[width=17em]{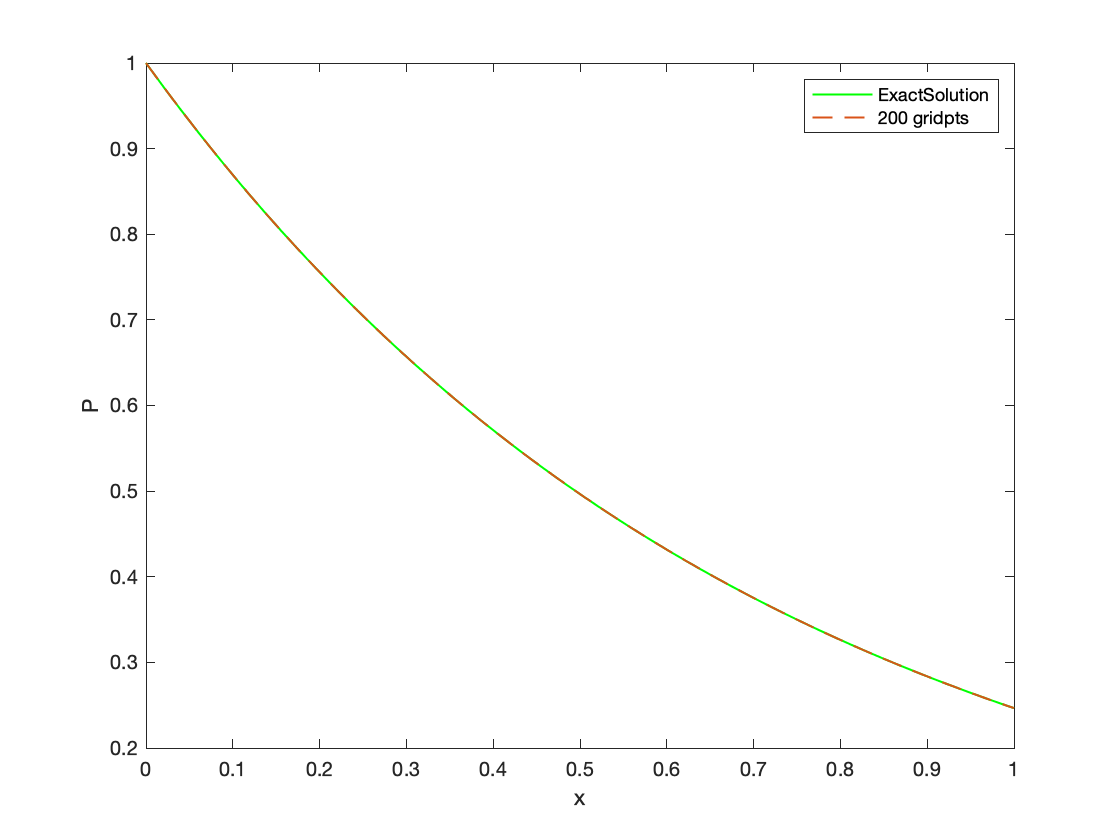} \\
    \small Moving Equilibrium: Energy &
      \small Moving Equilibrium: Pressure.
  \end{tabular}
  \vspace*{8pt}
  \caption{Solution of one-dimensional moving equilibrium. \label{fig.5.1-3}}
\end{figure}

In this experiment, we are interested in preserving the following moving equilibrium state with a nonlinear gravitational field $\phi(x) = \text{exp}(-\text{exp}(x)+\gamma(\text{exp}(-\gamma x)))$,
\begin{equation}
\begin{split}
    \rho(x) &= \rho_0 \text{exp}(-\frac{\rho_0g}{p_0}x),\\
    u(x) &= \text{exp}(x),\\
    p(x) &= \text{exp}(-\frac{\rho_0g}{p_0}x)^\gamma.
\end{split}
\end{equation}
with $\rho_0=1$ and $p_0=1$. The initial conditions are accordingly in the form: 
\begin{equation}
\begin{split}
    \rho(x) & = \text{exp}(-x), \\
    u(x) & = \text{exp}(x), \\
    p(x) & = \text{exp}(-\gamma x).
\end{split}
\end{equation}
The detailed introduction of the initial condition can be found in \cite{ref12}. The stationary solution used here is the equilibrium state itself and we apply outflow boundary conditions as in the first test case we presented earlier. We compute the solution in the interval $[0,1]$ on 200 grid points until the final time $t=10$ and compare it to the exact solution. The obtained results are shown in figure \ref{fig.5.1-3}. The displayed results show excellent matching with the exact analytical solution, thus endorsing the well-balanced characteristic of the developed numerical scheme.
\subsection{The Shock Tube Problem} \label{subsection5.1.4}

\begin{figure}[th]
  \centering
  \begin{tabular}{ c @{\quad} c }
    \includegraphics[width=18em]{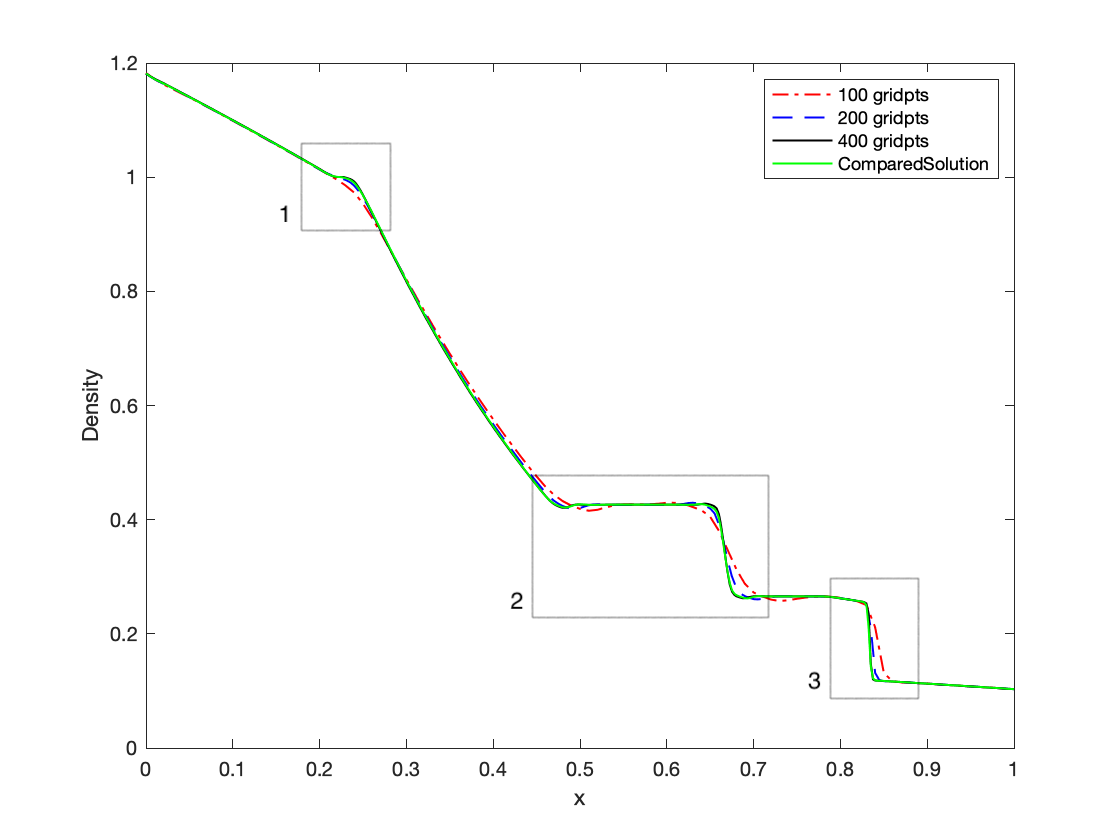} &
      \includegraphics[width=18em]{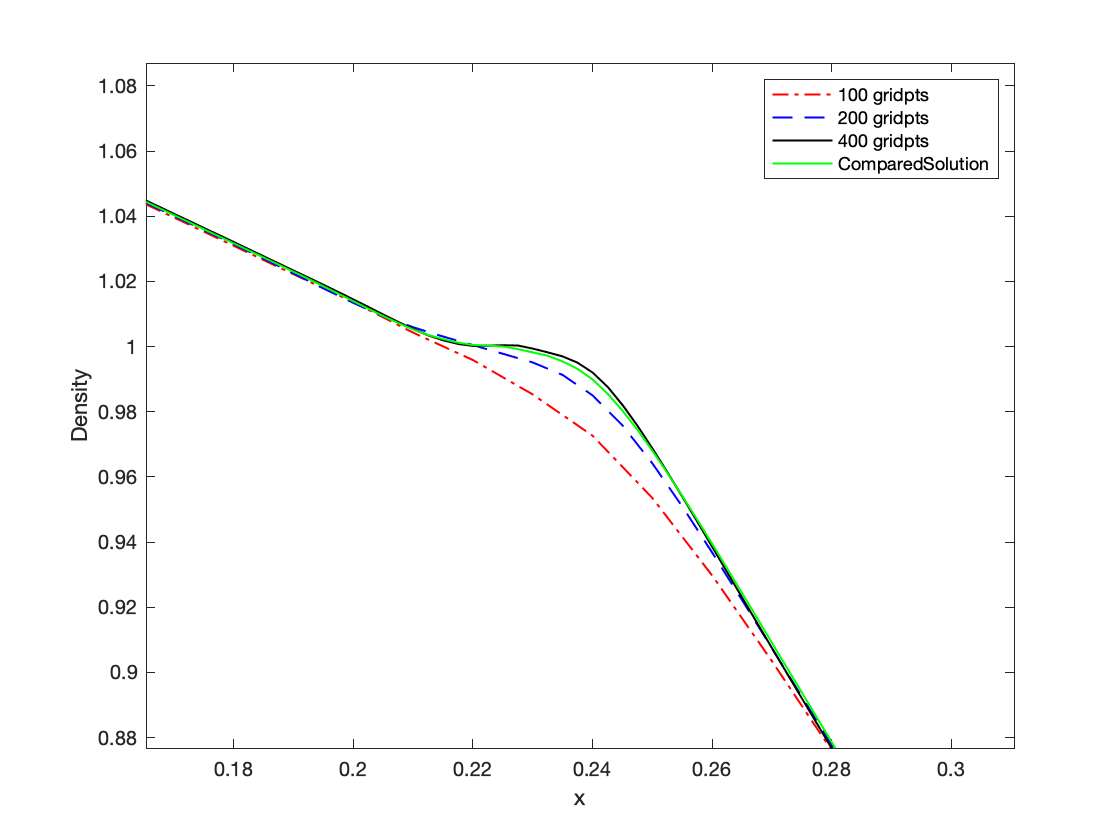} \\
    \small Density &
      \small Zoom in on the first block\\
     \includegraphics[width=18em]{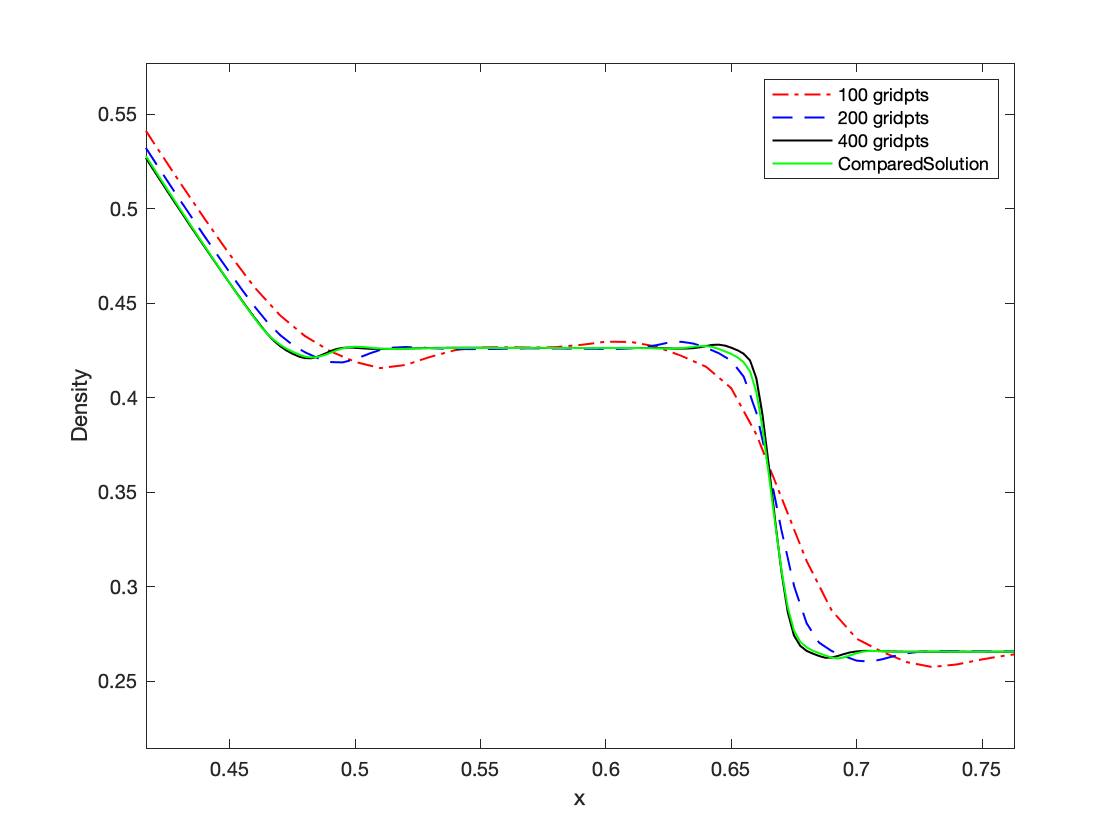} &
      \includegraphics[width=18em]{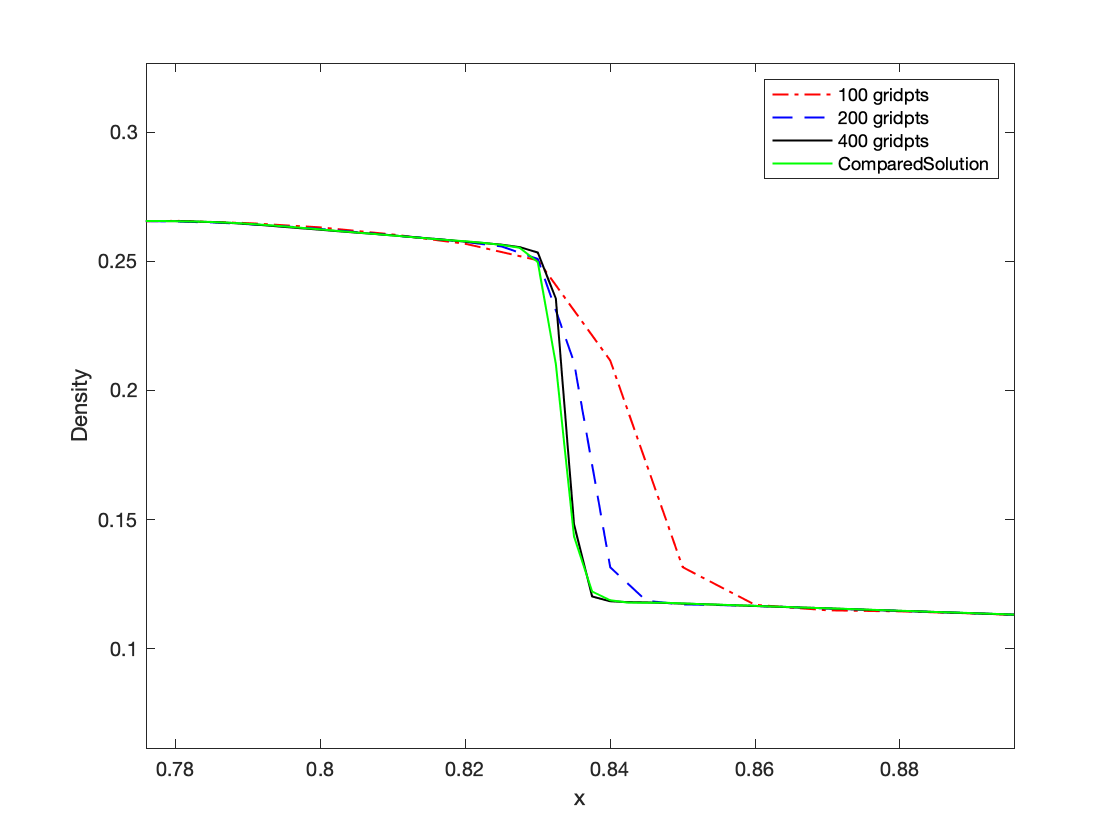} \\
    \small Zoom in on the second block&
      \small Zoom in on the third block 
  \end{tabular}
  \vspace*{8pt}
  \caption{Results of 1D shock tube problem. \label{fig.5-1.4.1}}
\end{figure}

\begin{figure}[th]
  \centering
  \begin{tabular}{ c @{\quad} c }
    \includegraphics[width=18em]{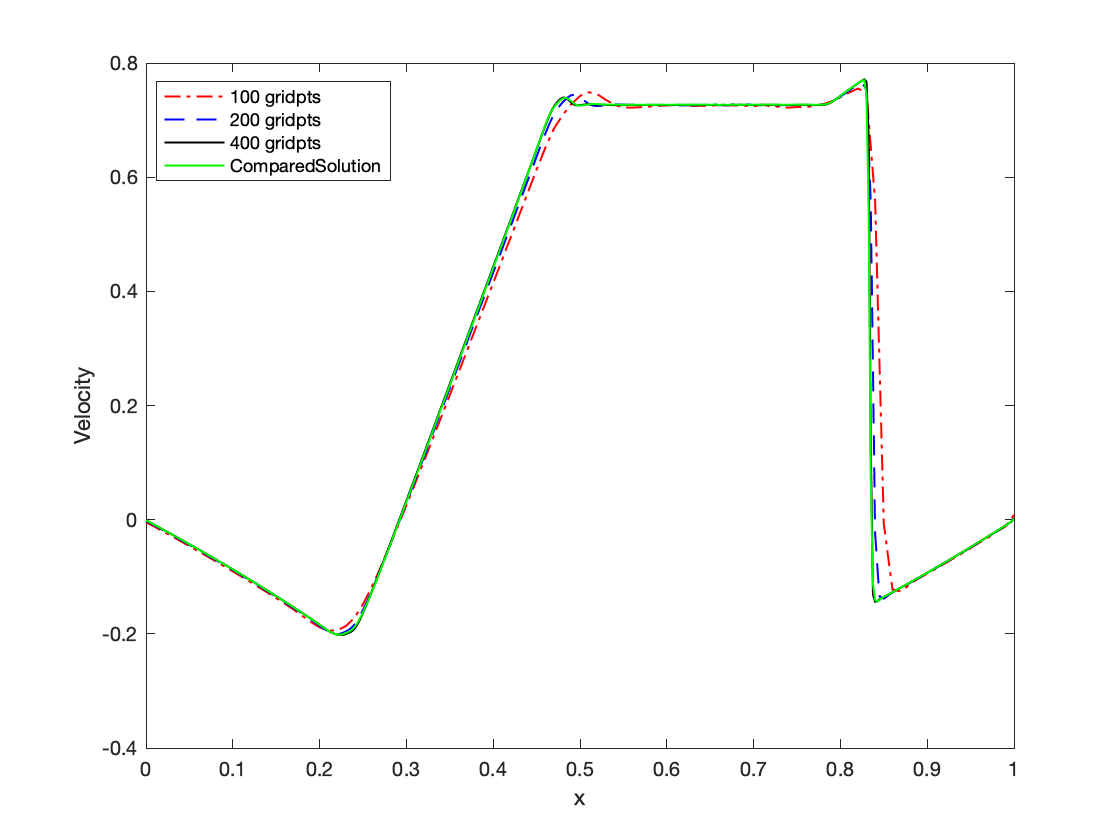}&
      \includegraphics[width=18em]{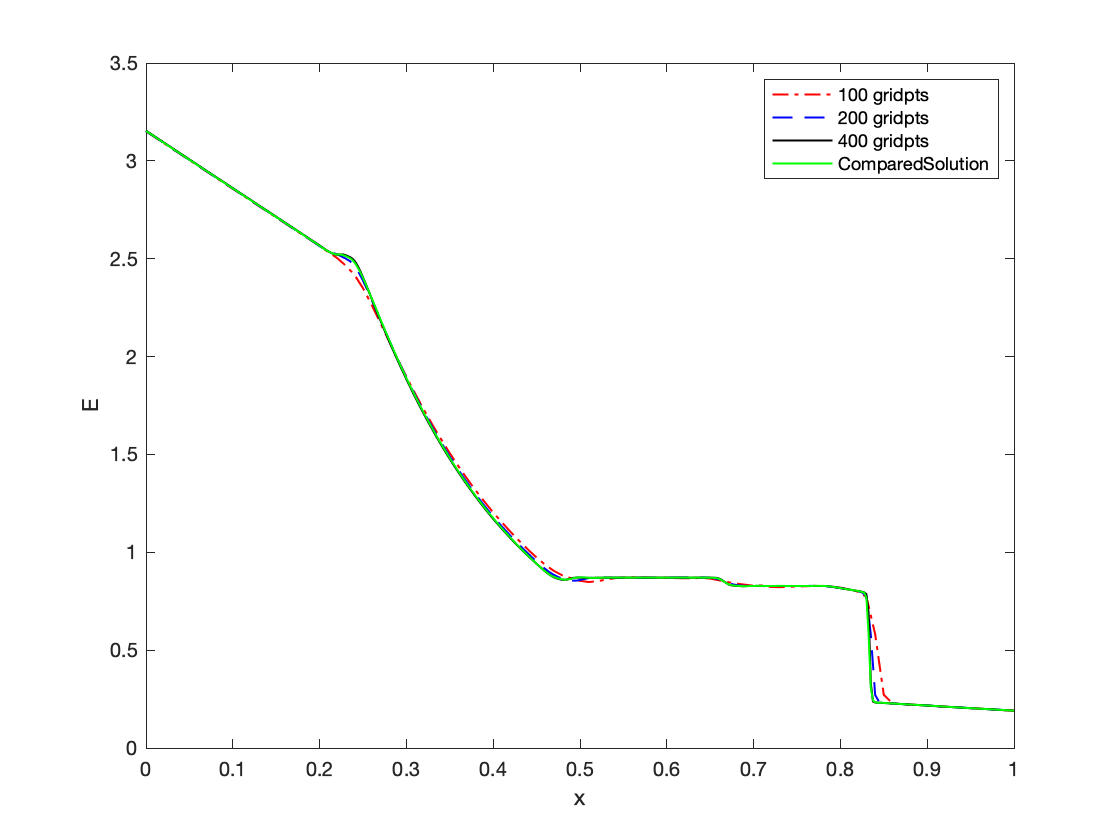} \\
    \small Velocity &
      \small Energy\\
     \includegraphics[width=18em]{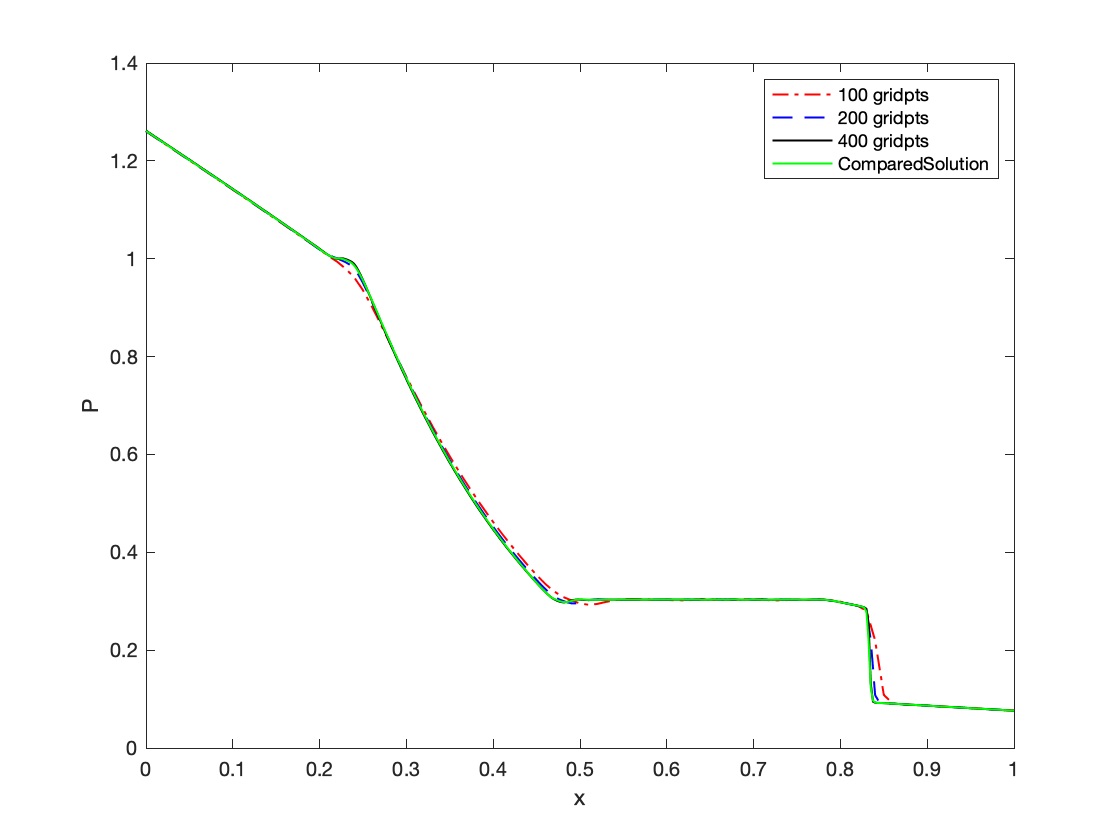}&
       \\
    \small Pressure&
      \small  
  \end{tabular}
  \vspace*{8pt}
  \caption{Results of 1D shock tube problem: velocity, energy and pressure. \label{fig.5-1.4.2}}
\end{figure}

In this example, we consider the shock tube problem with the definition from \cite{ref10}, and
we compare the obtained results with the solution from \cite{ref4}.
The initial state of the shock tube problem is given by
\begin{equation}
\begin{split}
    \rho(x) & =
    \left \{
    \begin{aligned}
        1,\qquad\;\;\; & \text{if}\;\; x\leq0.5, \\
        0.125,\quad & \text{otherwise},
    \end{aligned}
    \right. \\
    u(x) & = 0,\\
    p(x) & = 
    \left \{
    \begin{aligned}
        1,\quad\;\;\; & \text{if}\;\; x\leq0.5, \\
        0.1, \quad & \text{otherwise},
    \end{aligned}
    \right.
\end{split}
\end{equation}
and the gravitational field is defined with $\phi_x=g=1$. The stationary solution $\tilde q$ considered in this case is the isothermal equilibrium. \\
We use the computational interval $[0,1]$ with the reflecting boundary condition and we compute the solution on 100, 200, 400 grid points until the final time $t=0.2$. The reference solution is obtained from \cite{ref4} and computed on 400 grid points. In figure \ref{fig.5-1.4.1}, we show the results of the density and zoom in on the shocks. In figure \ref{fig.5-1.4.2}, we show the results of velocity, energy and pressure.
The reported results confirm the grid convergence and show a good comparison with the reference solution. 

\section{Conclusion}\label{sec:conc}
In this paper, we presented a new numerical scheme for the approximating solution of the in-homogeneous conservation laws. The numerical scheme we considered blends the central Kurganov-Tadmor scheme with the Deviation method.  Based on these two methods, we constructed a new second-order well-balanced scheme for the in-homogeneous Euler equations with gravity in one dimension. In addition, we showed that for the semi-discrete scheme applied to a scalar conservation law the TVD property holds in the regime of small deviations.

We applied our new scheme to the Euler equations with gravitational source term and successfully solved several classical problems. The solutions of the problem we obtained using our well-balanced numerical scheme are in excellent agreement with corresponding results from the literature. Still we noticed in a shock tube problem that there are some small oscillations in the vicinity of the shocks, especially when the grid is coarse. Increasing the number of grid points eliminates these oscillations. Extensions of the proposed scheme to the two-dimensional case are currently under investigations.

\section*{Acknowledgement}

The authors thank Praveen Chandrashekar (Tata Insitute, Centre for Applicable Mathematics, Bangalore, India) for very helpful discussions. 

\newcommand{\etalchar}[1]{$^{#1}$}

\appendix

\section*{Appendix} 

In this section, we provide the detailed calculations of the intermediate value $\overline{w}^{n+1}_{j+\frac{1}{2}}$ at evolution step in section 2.2.2. 
\par
Firstly, we compute the left-hand side ($LHS$) of equation \eqref{3.19},
\begin{align*}  
LHS=&\notag  \\
    \int^{x^n_{j+\frac{1}{2},r}}_{x^n_{j+\frac{1}{2},l}} & \bigg\{ \Big[f(\Delta q(x,t^n)+\tilde{q}(x,t^n))-f(\tilde{q}(x,t^n)) \Big]\,dt - \Delta q(x,t^n) \,dx \bigg\}  \notag  \\
    & +\int^{t_{n+1}}_{t_n} \bigg \{ \Big[f(\Delta q(x^n_{j+\frac{1}{2},r},t)+\tilde{q}(x^n_{j+\frac{1}{2},r},t))-f(\tilde{q}(x^n_{j+\frac{1}{2},r},t)) \Big]\,dt  \notag \\
    &\qquad\qquad\qquad\qquad\qquad\qquad\qquad\qquad\qquad\qquad\quad  - \Delta q(x^n_{j+\frac{1}{2},r},t)\, dx \bigg\}     \notag \\
    & +\int^{x^n_{j+\frac{1}{2},l}}_{x^n_{j+\frac{1}{2},r}} \bigg\{ \Big[f(\Delta q(x,t^{n+1})+\tilde{q}(x,t^{n+1}))-f(\tilde{q}(x,t^{n+1}))\Big]\,dt \notag \\
    &\qquad\qquad\qquad\qquad\qquad\qquad\qquad\qquad\qquad\qquad\quad   - \Delta q(x,t^{n+1}) \,dx \bigg\}  \notag  \\
    & +\int^{t_n}_{t_{n+1}} \bigg\{ \Big[f(\Delta q(x^n_{j+\frac{1}{2},l},t)+\tilde{q}(x^n_{j+\frac{1}{2},l},t))-f(\tilde{q}(x^n_{j+\frac{1}{2},l},t)) \Big]\,dt \notag \\
    &\qquad\qquad\qquad\qquad\qquad\qquad\qquad\qquad\qquad\qquad\quad   - \Delta q(x^n_{j+\frac{1}{2},l},t) \,dx \bigg\} \notag  
    \end{align*}
Rearranging equation \eqref{3.19}  and simplifying LHS we obtain:
\begin{equation}  \label{3.21}
\begin{split}
    \int^{x^n_{j+\frac{1}{2},r}}_{x^n_{j+\frac{1}{2},l}}  \Delta & q(x,t^{n+1}) dx  \\
    = & \int^{x^n_{j+\frac{1}{2},r}} _{x^n_{j+\frac{1}{2},l}} \Delta q(x,t^n) dx \\
    & - \int^{t_{n+1}}_{t_n}  [f(\Delta q(x^n_{j+\frac{1}{2},r},t)+\tilde{q} (x^n_{j+\frac{1}{2},r},t))-f(\tilde{q}(x^n_{j+\frac{1}{2},r},t))]dt \\
    & + \int^{t_{n+1}}_{t_n}  [f(\Delta q(x^n_{j+\frac{1}{2},l},t)+\tilde{q}(x^n_{j+\frac{1}{2},l},t))-f(\tilde{q}(x^n_{j+\frac{1}{2},l},t))]dt \\
    & + \int^{t^{n+1}}_{t^n} \int^{x^n_{j+\frac{1}{2},r}}_{x^n_{j+\frac{1}{2},l}} S(\Delta q,x) \;dxdt.
\end{split}
\end{equation}
Next, we define $\Delta x^n_{j+\frac{1}{2}}=x^n_{j+\frac{1}{2},r}-x^n_{j+\frac{1}{2},l} = 2a^n_{j+\frac{1}{2}}\Delta t$, and denote $f(\Delta q(x^n_{j+\frac{1}{2},r},t)+\tilde{q} (x^n_{j+\frac{1}{2},r},t))$ by $f((\Delta q+\tilde q)(x^n_{j+\frac{1}{2},r},t))$.\quad In order to get the average over  the interval $U^n_{j+\frac{1}{2}}$, which we denote by $\overline{w}^{n+1}_{j+\frac{1}{2}}$, we divide on both sides of equation \eqref{3.21} by the length  $\Delta x^n_{j+\frac{1}{2}}$ of the interval $U^n_{j+\frac{1}{2}}$, to obtain
\begin{equation}  \label{eqA.2}
\begin{split}
    \overline{w}^{n+1}_{j+\frac{1}{2}} 
    :=  & \frac{1}{\Delta x^n_{j+\frac{1}{2}}} \int^{x^n_{j+\frac{1}{2},r}}_{x^n_{j+\frac{1}{2},l}}  \Delta  q(x,t^{n+1}) dx  \\
    = & \frac{1}{\Delta x^n_{j+\frac{1}{2}}} \int^{x^n_{j+\frac{1}{2},r}} _{x^n_{j+\frac{1}{2},l}} \Delta q(x,t^n) dx \\
    & - \frac{1}{\Delta x^n_{j+\frac{1}{2}}} \int^{t_{n+1}}_{t_n} [f((\Delta q+\tilde q)(x^n_{j+\frac{1}{2},r},t)) - f(\tilde q (x^n_{j+\frac{1}{2},r},t))] \\
    & +\frac{1}{\Delta x^n_{j+\frac{1}{2}}} \int^{t_{n+1}}_{t_n} [f((\Delta q + \tilde q)(x^n_{j+\frac{1}{2},l},t)) - f(\tilde q (x^n_{j+\frac{1}{2},l},t))]dt \\
    & + \frac{1}{\Delta x^n_{j+\frac{1}{2}}} \int^{t^{n+1}}_{t^n} \int^{x^n_{j+\frac{1}{2},r}}_{x^n_{j+\frac{1}{2},l}} S(\Delta q,x) \;dxdt \\
    = & \frac{1}{\Delta x^n_{j+\frac{1}{2}}} \left [
    \int^{x^n_{j+\frac{1}{2},r}}_{x^n_{j+\frac{1}{2},l}} \Delta q(x,t^n) dx + F_U + S_U \right ],
\end{split}
\end{equation}
where   
\begin{equation} 
\begin{split}
    F_U
    = & - \int^{t_{n+1}}_{t_n} [f((\Delta q+\tilde{q}) (x^n_{j+\frac{1}{2},r},t))-f(\tilde{q}(x^n_{j+\frac{1}{2},r},t))]dt \\
    &+ \int^{t_{n+1}}_{t_n} [f(\Delta q+\tilde{q})(x^n_{j+\frac{1}{2},l},t))-f(\tilde{q}(x^n_{j+\frac{1}{2},l},t))]dt, 
\end{split}
\end{equation}
and
\begin{equation} 
    S_U = \int^{t^{n+1}}_{t^n} \int^{x^n_{j+\frac{1}{2},r}}_{x^n_{j+\frac{1}{2},l}} S(\Delta q,x) \;dxdt
\end{equation}

To evaluate the integral occurring in the last line of equation \eqref{eqA.2}, we apply the midpoint quadrature rule and we use the piecewise linearity of the reconstructions, thus we obtain
\begin{equation}  \label{3.27}
\begin{split}
    \frac{1}{\Delta x^n_{j+\frac{1}{2}}} & \int^{x^n_{j+\frac{1}{2},r}} _{x^n_{j+\frac{1}{2},l}} \Delta q(x,t^n) dx  \\
    & = \frac{1}{\Delta x^n_{j+\frac{1}{2}}} \Bigg [
    \int^{x^n_{j+\frac{1}{2}}}_{x^n_{j+\frac{1}{2},l}} \Delta q(x,t^n) dx
    +\int^{x^n_{j+\frac{1}{2},r}}_{x^n_{j+\frac{1}{2}}} \Delta q(x,t^n) dx \Bigg ]   \\
    & = \frac{1}{\Delta x^n_{j+\frac{1}{2}}} 
    \Bigg [ \frac{\Delta x^n_{j+\frac{1}{2}}}{2} Q_j(x^n_{j+\frac{1}{2},lm},t^n)
    + \frac{\Delta x^n_{j+\frac{1}{2}}}{2} Q_{j+1}(x^n_{j+\frac{1}{2},rm},t^n) \Bigg ]  \\
    & = \frac{1}{\Delta x^n_{j+\frac{1}{2}}} 
    \Bigg [ \frac{\Delta x^n_{j+\frac{1}{2}}}{2} 
    \bigg ( (\Delta q)^n_j + (x^n_{j+\frac{1}{2},lm} - x_j ) ((\Delta q)_x)^n_j)    \\
    & \qquad \qquad \qquad \qquad   + (\Delta q)^n_{j+1} + (x^n_{j+\frac{1}{2},rm} - x_{j+1} ) ((\Delta q)_x)^n_{j+1}) \bigg ) 
     \Bigg ]      \\
    & = \frac{1}{\Delta x^n_{j+\frac{1}{2}}}
    \Bigg [ \frac{\Delta x^n_{j+\frac{1}{2}}}{2} 
    \bigg( (\Delta q)^n_j + (\frac{\Delta x}{2} - \frac{\Delta x^n_{j+\frac{1}{2}}}{4}) ((\Delta q)_x)^n_j)       \\
    & \qquad \qquad \qquad \qquad   + (\Delta q)^n_{j+1} - (\frac{\Delta x}{2} - \frac{\Delta x^n_{j+\frac{1}{2}}}{4}) ((\Delta q)_x)^n_{j+1}) \bigg )   \Bigg ]  \\
    & = \frac{(\Delta q)^n_j + (\Delta q)^n_{j+1}}{2}
    + \frac{\Delta x - \frac{\Delta x^n_{j+\frac{1}{2}}}{2}}{4} \bigg[ ((\Delta q)_x)^n_j - ((\Delta q)_x)^n_{j+1} \bigg],
\end{split}
\end{equation}
where $x^n_{j+\frac{1}{2},lm}$ denotes the midpoint of the interval $[x^n_{j+\frac{1}{2},l}, x^n_{j+\frac{1}{2}}]$; similarly $x^n_{j+\frac{1}{2},rm}$ denotes the midpoint of the interval $[x^n_{j+\frac{1}{2}}x^n_{j+\frac{1}{2},r}]$.


\end{document}